\documentclass[aip,graphicx]{revtex4-1}

\draft 

\begin{document}


\title{Reflection identities of harmonic sums of  weight four } 



\author{Alex Prygarin}
\email[]{alexanderp@ariel.ac.il}
\thanks{This paper is dedicated to memory of Lev Lipatov}
\affiliation{Department of Physics
\\
 Ariel University\\
Ariel, 40700, Israel}



\begin{abstract}
We consider the reflection identities for harmonic sums at weight four. We decompose a product of two harmonic sums with mixed pole structure into a linear combination of terms each having a pole at either negative or positive values of the argument. The pole decomposition demonstrates how the product of two simpler harmonic sums can build more complicated harmonic sums at higher weight. We list a minimal irreducible bilinear set of reflection identities at weight four which present the main result of the paper. We also discuss how other trilinear and  quartic reflection identities can be easily constructed from our result with the use of well known shuffle relations for harmonic sums.
\end{abstract}


\maketitle 


\section{Introduction}

In this paper we continue discussion of our previous study~\cite{Prygarin:2018tng} regarding the reflection identities of harmonic sums, where a product of two harmonic sums of argument $z$ and $-1-z$ is expressed through a linear combination of other harmonic sums of the same arguments, i.e.
\begin{eqnarray}
S_{a_1,a_2,...}(z)S_{b_1, b_2, ...}(-1-z)=S_{c_1,c_2,...}(z)+...+S_{d_1,d_2,...}(-1-z)+... 
\end{eqnarray}

The reflection identities at weight two presented here are not new and were known long time ago in the context of functions related to the Euler Gamma function. To the best of our knowledge they appear the  eraliest in  Chapter 20 of the book by Nielsen~\cite{handbuch} and then were related to the harmonic sums~(see eqs.6.11-6.15 of the paper by  Blumlein~\cite{Blumlein:2009ta}). At weight three they were recently calculated by the author~\cite{Prygarin:2018tng}. This paper deals with weight four. 

 The harmonic sums have pole singularities at negative integers. The reflection identities present a pole separation for a product of two sums with mixed pole structure. 
 We call those functional relations the reflection identities because the argument of the harmonic sums is reflected with respect to the point $\frac{z+(-1-z)}{2}=-\frac{1}{2}$. The reflection identities up to  weight three were published in our previous study~\cite{Prygarin:2018tng} and here we present them at weight four. 
 
The harmonic sums are defined through a nested summation with their argument being the upper limit in the outermost sum~\cite{HS1,Vermaseren:1998uu,Blumlein:1998if,Remiddi:1999ew}
\begin{eqnarray}\label{defS}
S_{a_1,a_2,...,a_k}(n)=  \sum_{n \geq i_1 \geq i_2 \geq ... \geq i_k \geq 1 }   \frac{\mathtt{sign}(a_1)^{i_1}}{i_1^{|a_1|}}... \frac{\mathtt{sign}(a_k)^{i_k}}{i_k^{|a_k|}}
\end{eqnarray}
In this paper we consider the harmonic sums with only real integer values of $a_i$, which build the alphabet of the possible negative and positive indices.  
In Eq.~(\ref{defS}) $k$ is  the depth and $w=\sum_{i=1}^{k}|a_i|$ is the weight of the harmonic sum $S_{a_1,a_2,...,a_k}(n)$. 

The indices of harmonic sums $a_1,a_2,...,a_k$ can be either positive or negative integers and label uniquely $S_{a_1,a_2,...,a_k}(n)$ for any given 
weight. However there is no unique way of building the functional basis for a given weight because the harmonic sums are subject to so called shuffle relations, where a linear combination of $S_{a_1,a_2,...,a_k}(n)$ with the same argument but all possible permutations of indices can be expressed through 
a non-linear combinations of harmonic sums at lower weight.
 There is also some freedom in choosing the irreducible minimal set of $S_{a_1,a_2,...,a_k}(n)$ that builds those non-linear combinations.  
The shuffle relations make a connection between the linear and non-linear combinations of the harmonic sums of the same argument. For example,  the shuffle relation at depth two reads
\begin{eqnarray}\label{shuffle}
S_{a,b}(z)+S_{b,a}(z)= S_{a}(z) S_{b}(z)+S_{\textrm{ sign}(a) \textrm{ sign}(b)(|a|+|b|)}(z)
\end{eqnarray}
The shuffle relations of the harmonic sums is closely connected to the shuffle algebra of the  harmonic polylogarithms~\cite{Remiddi:1999ew}.

There is another type of identity called the duplication identities where a combination of harmonic sums of argument $n$ can be expressed through 
a harmonic sum of the argument $2 n$. The duplication identities introduce another freedom in choosing the functional basis.

  In this paper we consider the analytic continuation of the harmonic sums to from positive integer values of the argument to the complex plane denoted by $\bar{S}_{a_1,a_2,...}^{+}(z)$~(this notation was introduced by Kotikov and Velizhanin~\cite{Kotikov:2005gr}). The analytic continuation is done in terms of the Mellin transform of corresponding Harmonic Polylogarithms and was recently used by Gromov, Levkovich-Maslyuk and Sizov~\cite{Gromov:2015vua,Alfimov:2018cms} and Caron Huot and Herraren~\cite{Caron-Huot:2016tzz} for expressing the  eigenvalue Balitsky-Fadin-Kuraev-Lipatov~(BFKL)  equation using the principle of Maximal Transcedentality~\cite{Kotikov:2006ts} in super Yang-Mills $\mathcal{N}=4$ field theory. We plan to use their results together with analysis done by one of the authors and collaborators~\cite{Bondarenko:2015tba,Bondarenko:2016tws} to understand the general structure of the   BFKL equation in QCD and beyond.

 The Mellin transform allows to make the analytic continuation to the complex plane. For example, consider the harmonic sum 
 \begin{eqnarray}\label{s1}
 S_{-1}(z)=\sum_{k=1}^{z} \frac{(-1)^k}{k}
\end{eqnarray}  
The corresponding Mellin transform reads
\begin{eqnarray}
\int_0^1 \frac{1}{1+x} x^z= (-1)^z \left(S_{-1}(z)+\ln 2\right)
\end{eqnarray}
One can see that $S_{-1}(z)$ on its own is not an analytic function because of the term $(-1)^z $ and we impose that we start from even integer values of the argument $z$. In this case we define its analytic continuation from even positive integers to all positive integers through 
\begin{eqnarray}
\bar{S}_{a_1,a_2,...}^{+}(z)=(-1)^z S_{-1}(z)+((-1)^z-1) \ln 2
\end{eqnarray} 
 and thus we can write 
 \begin{eqnarray}
 \bar{S}_{a_1,a_2,...}^{+}(z)=\int_0^1 \frac{1}{1+x} x^z-\ln 2
\end{eqnarray}
This way we defined  $\bar{S}_{a_1,a_2,...}^{+}(z)$ using the Mellin transform of ratio function $ \frac{1}{1+x}$. In more complicated cases of other harmonic sums one includes also Harmonic Polylogarithms on top of the ratio functions, but the general procedure is very similar and largely covered in  a number of publications~\cite{AblingerThesis,Ablinger:2011te,Blumlein:2009ta,Blumlein:1998if,Blumlein:2009fz}.

 It is worth mentioning that there is another analytic continuation for the harmonic sum, from odd positive integer values of the argument, which is different for harmonic sums with at least one negative index and  denoted by  $\bar{S}_{a_1,a_2,...}^{-}(z)$. 
Both analytic continuations are equally valid. 
 Our goal is to find a closed expression of the BFKL eigenvalue for all possible values of anomalous dimension and conformal spin, so that  we follow the notation of Gromov, Levkovich-Maslyuk and Sizov~\cite{Gromov:2015vua}, and use $\bar{S}^{+}(z)$ throughout the text. For simplicity of presentation in this paper we write everywhere $S_{a_1,a_2,...}(z)$ instead of $\bar{S}_{a_1,a_2,...}^{+}(z)$.   

 As it was already mentioned there is no unique way in defining a minimal irreducible set of harmonic sums due to the functional relations between them. 
 For example, one can use the shuffle relations and them the minimal irreducible basis would include quadratic terms $S_{a_1}(z)S_{a_2}(z)$ in place of either $S_{a_1,a_2}(z)$ or $S_{a_2,a_1}(z)$. It is convenient to use shuffle relations to remove from the minimal basis the harmonic sums with the first index being equal $1$, because those are divergent as $z \to \infty$. Then, the remaining harmonic sums give transcendental constants at $z \to \infty$. Most of constants are reducible and one is free to choose an irreducible set of transcendental  constants at any given weight. We use the set implemented in the HarmonicSums package. The irreducible constants are given by 
\begin{eqnarray}\label{C1}
C_1=\{\log (2)\}
\end{eqnarray}
and 
\begin{eqnarray}\label{C2}
C_2=\left\{\pi ^2,\log ^2(2)\right\}
\end{eqnarray}
and 
\begin{eqnarray}\label{C3}
C_3=\left\{\pi ^2 \log (2),\log ^3(2),\zeta _3\right\}
\end{eqnarray}
as well as 
\begin{eqnarray}\label{C4}
C_4=\left\{\pi ^4,\pi ^2 \log ^2(2),\log ^4(2),\text{Li}_4\left(\frac{1}{2}\right),\zeta _3 \log
   (2)\right\}
\end{eqnarray}
 where $C_w$ stands for a minimal set of irreducible constants at given weight $w$. There is only one of those  at $w=1$, two at $w=2$, three at $w=3$ and five irreducible constants at weight $w=4$.

We choose to use a linear minimal set of the harmonic sums to represent our results. In this set we do not apply  shuffle relations and thus all the terms of the basis are linear in $S_{a_1,a_2,...}(z)$. This choice is dictated mostly by a convenience and was also used by Caron Huot and Herraren~\cite{Caron-Huot:2016tzz} on which we would like to rely   in  our future calculations. The minimal linear set of harmonic sums we use is  as follows
\begin{eqnarray}\label{B_1}
B_1=\left\{S_{-1},S_1\right\}
\end{eqnarray}
and 
\begin{eqnarray}\label{B_2}
B_2&=&\left\{S_{-2},S_2,S_{-1,1},S_{1,-1},S_{1,1},S_{-1,-1}\right\}
\end{eqnarray}
and 
\begin{eqnarray}\label{B_3}
B_3&=&\left\{S_{-3},S_3,S_{-2,-1},S_{-2,1},S_{2,-1},S_{2,1},S_{-1,1,-1},S_{-1,1,1}, S_{1,-2},S_{1,2},S_{1
   ,-1,-1},S_{1,-1,1},\right.
\nonumber
\\
&& 
\left.
S_{1,1,-1},S_{1,1,1},S_{-1,-2},S_{-1,2},S_{-1,-1,-1},S_{-1,-1,1}\right\}
\end{eqnarray}
as well as
\begin{eqnarray}\label{B_4}
 B_4&=&\left\{S_{-4},S_4,S_{-3,-1},S_{-3,1},S_{2,-2},S_{3,-1},S_{3,1},S_{-2,-1,-1},S_{-2,-1,1},S_{-2,1,-
   1},S_{-2,1,1},
   \right.
\nonumber
\\
&& 
\left.
S_{2,-1,-1},S_{2,-1,1},S_{2,1,-1},S_{2,1,1},S_{-1,1,-1,-1},S_{-1,1,-1,1},S_{-1,1
   ,1,1},S_{1,-3},S_{1,3},S_{1,-2,-1},
   \right.
\nonumber
\\
&& 
\left.
S_{1,-2,1},S_{1,-1,-2},S_{1,-1,2},S_{1,1,-2},S_{1,1,2},S_{1
   ,2,-1},S_{1,2,1},S_{1,-1,-1,-1},S_{1,-1,-1,1},S_{1,-1,1,-1},
   \right.
\nonumber
\\
&& 
\left.
S_{1,-1,1,1},S_{1,1,-1,-1},S_{1,1,
   -1,1},S_{1,1,1,-1},S_{1,1,1,1},S_{-1,-3},S_{-1,3},S_{-1,-2,-1},S_{-1,-2,1},
   \right.
\nonumber
\\
&& 
\left.
S_{-1,-1,-2},
S_{-1,
   -1,2},S_{-1,1,-2},S_{-1,1,2},S_{-1,2,-1},S_{-1,2,1},S_{-1,-1,-1,-1},S_{-1,-1,-1,1},
   \right.
\nonumber
\\
&& 
\left.
S_{-1,-1,1,
   -1},
S_{-1,-1,1,1},
S_{-1,1,1,-1},S_{-2,-2},S_{-2,2},S_{2,2}\right\}
\end{eqnarray}

A comprehensive discussion on harmonic sums,  irreducible constants, functional identities and possible choice of the minimal set of functions at given weight is presented by J.~Ablinger~\cite{AblingerThesis}. 
 In this paper we focus only at the reflection identities for harmonic sums at weight $w=4$ analytically continued from even positive points to complex plane.  In the next Section we discuss them in more details along the method we use in our calculations.

\section{Reflection identities}

The reflection identities at weight $w=4$ are obtained by taking a product of harmonic sums of argument $z$ at weight $w=1$ and harmonic sums of 
argument $-1-z$
at weight $w=3$, i.e. $B_1 \otimes \bar{B}_3$, and also by taking a product  of harmonic sums of argument $z$ and $-1-z$ at weight $w=2$, 
i.e. $B_2 \otimes \bar{B}_2$.

The number of basis harmonic sums in $B_1$, $B_2$ and $B_3$ is given by 
\begin{eqnarray}\label{L1L2L3}
\text{Length}(B_1)=2, \;\;\; \text{Length}(B_2)=6,  \;\;\;  \text{Length}(B_3)=18
\end{eqnarray}
so that the number of elements in the products $B_1 \otimes \bar{B}_3$ and  $B_2 \otimes \bar{B}_2$ reads 
\begin{eqnarray}
B_1 \otimes \bar{B}_3=2 \times 18=36
\end{eqnarray}
and 
\begin{eqnarray}
B_2 \otimes \bar{B}_2=\frac{6 \times  (6-1)}{2}+6=21,
\end{eqnarray}
resulting in the total number of irreducible reflections identities at weight $w=4$ being equal to $21+36=57$.

In order to calculate the reflection identities at weight $w=4$ we use the  basis harmonic sums at $w=4$ listed in Eq.~(\ref{B_4}) together with basis 
 harmonic sums at lower weight listed in the work of J.~Ablinger~(\ref{B_1})-(\ref{B_3}) multiplied by irreducible constants at corresponding weight listed in Eqs.~(\ref{C1})-(\ref{C3}).  This should be supplemented by irreducible constants at weight $w=4$ listed in Eqs.~(\ref{C4}).
The number of basis sums in $B_4$ equals
\begin{eqnarray}\label{L4}
 \text{Length}(B_4)=54
\end{eqnarray}
so that   the total number of terms in the expansion ansatz at $w=4$ 
\begin{eqnarray}
 B_4 + B_3 \otimes C_1+B_2 \otimes C_2+B_1 \otimes C_3+C_4 
\end{eqnarray}
 is given by
 \begin{eqnarray}\label{Lanz4}
\text{Length}(\text{ANZ}_4)=54 + 18 \times 1+6 \times 2+2  \times 3 + 5=95 
\end{eqnarray}
The full expansion ansatz at $w=4$ is given by 
\begin{eqnarray}\label{anz4}
\text{ANZ}_4&=&\left\{\pi ^4,\pi ^2 \log ^2(2),\log ^4(2),\text{Li}_4\left(\frac{1}{2}\right),S_{-4},S_{-3} \log
   (2),\pi ^2 S_{-2},S_{-2} \log ^2(2),
   \right.
\nonumber
\\
&& 
\left.
\pi ^2 S_{-1} \log (2),S_{-1} \log ^3(2),\pi ^2 S_1 \log
   (2),S_1 \log ^3(2),\pi ^2 S_2,S_2 \log ^2(2),
   \right.
\nonumber
\\
&& 
\left.
S_3 \log
   (2),S_4,S_{-3,-1},S_{-3,1},S_{-2,-2},\log (2) S_{-2,-1},\log (2)
   S_{-2,1},S_{-2,2},
   \right.
\nonumber
\\
&& 
\left.
S_{-1,-3},\log (2) S_{-1,-2},\pi ^2 S_{-1,-1},\log ^2(2) S_{-1,-1},\pi ^2
   S_{-1,1},\log ^2(2) S_{-1,1},
   \right.
\nonumber
\\
&& 
\left.
\log (2) S_{-1,2},S_{-1,3},S_{1,-3},\log (2) S_{1,-2},\pi ^2
   S_{1,-1},\log ^2(2) S_{1,-1},\pi ^2 S_{1,1},
   \right.
\nonumber
\\
&& 
\left.
\log ^2(2) S_{1,1},\log (2)
   S_{1,2},S_{1,3},S_{2,-2},\log (2) S_{2,-1},\log (2)
   S_{2,1},S_{2,2},S_{3,-1},
   \right.
\nonumber
\\
&& 
\left.
S_{3,1},
   S_{-2,-1,-1},S_{-2,-1,1},S_{-2,1,-1},S_{-2,1,1},S_{-1,-2,-1},
   S_{-1,-2,1},S_{-1,-1,-2},
   \right.
\nonumber
\\
&& 
\left.
\log (2) S_{-1,-1,-1},\log (2)
   S_{-1,-1,1},S_{-1,-1,2},S_{-1,1,-2},\log (2) S_{-1,1,-1},
   \right.
\nonumber
\\
&& 
\left.
\log (2)
   S_{-1,1,1},S_{-1,1,2},S_{-1,2,-1},S_{-1,2,1},S_{1,-2,-1},S_{1,-2,1},S_{1,-1,-2},
   \right.
\nonumber
\\
&& 
\left.
\log (2)
   S_{1,-1,-1},\log (2) S_{1,-1,1},S_{1,-1,2},S_{1,1,-2},\log (2) S_{1,1,-1},\log (2)
   S_{1,1,1}
   ,
   \right.
\nonumber
\\
&& 
\left.
S_{1,1,2},S_{1,2,-1},S_{1,2,1},S_{2,-1,-1},S_{2,-1,1},S_{2,1,-1},S_{2,1,1},S_{-1,-1,
   -1,-1},S_{-1,-1,-1,1},
   \right.
\nonumber
\\
&& 
\left.
S_{-1,-1,1,-1},S_{-1,-1,1,1},S_{-1,1,-1,-1},S_{-1,1,-1,1},S_{-1,1,1,-1},
   S_{-1,1,1,1},S_{1,-1,-1,-1},
   \right.
\nonumber
\\
&& 
\left.
S_{1,-1,-1,1},S_{1,-1,1,-1},S_{1,-1,1,1},S_{1,1,-1,-1},S_{1,1,-1,1
   },S_{1,1,1,-1},S_{1,1,1,1},
   \zeta _3 \log (2),
   \right.
\nonumber
\\
&& 
\left.
\zeta _3 S_{-1},\zeta _3 S_1\right\}
\end{eqnarray}
The expansion of the product of two functions of argument $z$ and argument $-1-z$ we search in terms of two sets of $\text{ANZ}_4$, one of argument $z$ and another one of argument $-1-z$. The total number of elements in this expression equals $95 \times 2-5=185$, where we remove redundant five 
constants at $w=4$ because they are the same for both arguments. We fix the $185$ free coefficients using pole expansion of the product
 $s_{a_1,a_2,..}(z)s_{b_1,b_2,..}(-1-z)$ around negative integers of $z=-5,...,-1$ and expanding to the second order of the expansion parameter. It turns out
 that to fix all $185$ free coefficients we need only expansion up to first order and we use the second order of the expansion to double check our results.
 We checked our results listed in the Appendix by a direct numerical calculation at the accuracy  $10^{-10}$. 

Below we give two examples of reflection identities, $ S_1(z)S_{2,1}(-1-z)$ for harmonic sums with positive indices and 
 $S_{-1}(z)    S_{-2,-1}(-1-z)$ for harmonic sums with negative indices.  All $57$  irreducible 
  reflection identities at weight four are listed in the Appendix. The two examples are 
 \begin{eqnarray}\label{S1S21}
     S_1(z)S_{2,1}(-1-z)& = & \frac{6 \zeta _2^2}{5}-S_2(z) \zeta _2+S_2(-1-z) \zeta _2+2 S_1(z) \zeta _3
     -2 \zeta
   _3 S_1(-1-z)
   +S_{3,1}(z) \nonumber
    \\
     &&-S_{3,1}(-1-z)-S_{2,1,1}(z)+S_{1,2,1}(-1-z)+S_{2,1,1}(-1-z)
   \end{eqnarray}
and 
 \begin{eqnarray}\label{Sm1sSm2m1}
    S_{-1}(z)    S_{-2,-1}(-1-z)& = & -\frac{\ln^4 ( 2)}{6}-S_{-2}(z) \ln^2(2)+S_2(z) \ln^2(2)+\zeta_2 \ln^2(2)
    \nonumber \\
   &&
   -S_{-2}(-1-z) \ln^2(2)
    +S_2(-1-z) \ln^2(2)-S_{-3}(z) \ln(2)
    \nonumber \\
   &&
   +S_3(z) \ln (2)    -6 \zeta _3 \ln (2)
   -S_{-3}(-1-z) \ln (2)
   +S_3(-1-z) \ln (2)
   \nonumber \\
   &&
   -2 S_{-2,-1}(z) \ln (2)-2 S_{-2,-1}(-1-z) \ln (2) +\frac{33\zeta _2^2}{20}
   -4 \text{Li}_4 \left(\frac{1}{2}\right)
   \nonumber \\
   &&
   +\frac{1}{2} S_{-2}(z) \zeta _2
      +\frac{1}{4} S_{-1}(z) \zeta
   _3
   +\frac{1}{2} \zeta _2 S_{-2}(-1-z)+\frac{1}{4} \zeta _3
   S_{-1}(-1-z)
   \nonumber \\
   &&
   +S_{3,-1}(z)+S_{3,-1}(-1-z)-S_{-2,-1,-1}(z)
   -S_{-2,-1,-1}(-1-z)
   \nonumber \\
   &&
   -S_{-1,-2,-1}(-1-z)
    \end{eqnarray} 

One can see that  the reflection identities for harmonic sums with negative indices are  more complicated than those with only positive indices and this happens mostly due to appearance of constant $\ln ( 2)$, which originates from sign alternating summation in $S_{-1}(z)$   absent for positive indices.

In the present paper we consider only bilinear reflection identities expressing a product of two harmonic sums of argument $z$ and $-1-z$ in terms of a linear combination of other sums of the same arguments. One can consider also trilinear and quartic identities, but whose reducible and  form   a linear combination of the bilinear identities presented in this paper. For example, we can consider a trilinear term
$S_1(z)S_{1}(-1-z)S_{2}(-1-z)$ and write it as 
\begin{eqnarray}\label{triples1sc1sc2}
S_1(z)S_{1}(-1-z)S_{2}(-1-z)=S_1(z)S_{1,2}(-1-z)+S_1(z)S_{2,1}(-1-z)-S_1(z)S_{3}(-1-z) \hspace{1cm}
\end{eqnarray}
where we used a shuffle identity from Eq.~(\ref{shuffle})
\begin{eqnarray}
S_{1,2}(z)+S_{2,1}(z)-S_{3}(z)=S_{1}(z)S_{2}(z)
\end{eqnarray}
The expression for $S_1(z)S_{1,2}(-1-z)$ is given in Eq.~(\ref{s1s12}) and for $S_1(z)S_{2}(-1-z)$  in Eq.~(\ref{s1s3}). Plugging those together with  Eq.~(\ref{S1S21}) into Eq.~(\ref{triples1sc1sc2}) we get 
\begin{eqnarray}
S_1(z)S_{1}(-1-z)S_{2}(-1-z)&=& -\zeta _2 S_{1,1}(-1-z)-\bar{S}_{1,3}-S_{2,2}(-1-z)-S_{3,1}(-1-z)
\nonumber
\\
&&
+2
   S_{1,1,2}(-1-z)+S_{1,2,1}(-1-z)+S_{2,1,1}(-1-z)
   \nonumber
\\
&&
+2 \zeta _2 S_2(-1-z)+\zeta _3 S_1(-1-z)+\zeta
   _2 S_{1,1}(z)+S_{3,1}(z)
   \nonumber
\\
&&
-S_{1,2,1}(z)-S_{2,1,1}(z)+\frac{4 \zeta _2^2}{5}-\zeta _2 S_2(z)+2 \zeta _3 S_1(z)
\end{eqnarray}
In a similar way one can can build any trilinear or quartic reflection identity using shuffle relations for harmonic sums and the bilinear reflection identities listed in the Appendix of this paper. All possible shuffle relations 
 required for the present discussion  are available in the HarmonicSums package by J.~Ablinger.  Shuffle relation before and after analytic continuation of the harmonic sums to the complex plane are the same.

\section{Conclusions}
We discuss the reflection identities for harmonic sums of weight four. There are $57$ irreducible bilinear identities listed in the Appendix. All other bilinear reflection identities are easily obtained by a trivial change of argument $z \leftrightarrow -1-z$. The trilinear and quartic identities for a product of three and four harmonic sums are obtained from the identities listed in the Appendix using the shuffle relations for harmonic sums. In our analysis we use the linear basis for harmonic sums and limit ourselves to harmonic sums analytically continued from even integer values of the argument to the complex plane. The analytic continuation from odd integers is beyond  the scope of the present study. 

In deriving the reflection identities presented in this paper we used Harmonic Sums package by J.~Ablinger~\cite{AblingerThesis}, HPL package by D.~ Maitre~\cite{Maitre:2005uu} and dedicated Mathematica package for pomeron NNLO eigenvalue by N.~ Gromov, F.~Levkovich-Maslyuk and G.~Sizov~\cite{Gromov:2015vua}. 

We expanded around positive and negative integer points  the  product of two harmonic sums $S_{a_1,a_2,..}(z)S_{b_1,b_2,...}(-1-z)$ and the functional basis built of pure Harmonic Sums with constants of relevant weight  listed in Ref.~\cite{AblingerThesis}. Then we compared the coefficients of the irreducible constants of a given weight and solved the resulting set of coefficient equations. We used higher order expansion to cross check our results.   The bilinear reflection identities presented here are derived from the pole expansion based on the Mellin transform and then checked them  against the shuffle identities and numerical calculations of the corresponding harmonic sums on the complex plane. 

We attach a Mathematica notebook with our results.

\section{Acknowledgements}
We would like to thank Fedor Levkovich-Maslyuk and Mikhail Alfimov for fruitful discussions on details of their calculations of NNLO BFKL eigenvalue~\cite{Gromov:2015vua,Alfimov:2018cms}. We are grateful to Simon Caron-Huot for explaining us the structure of his result on NNLO BFKL eigenvalue and his calculation techniques~\cite{Caron-Huot:2016tzz}.

We are indebted to Jochen Bartels for his hospitality and enlightening discussions during our stay at University of Hamburg where this project was initiated.

\newpage

\subsection{Appendix}
 Here we list below the irreducible reflection identities at weight $w=4$. In all our expressions we used the linear minimal set of harmonic sums given in Eqs.~(\ref{B_1})-(\ref{B_4}).

We use a compact notation where  $s_{a_1, a_2, ...}$ stands for $S_{a_1, a_2, ...}(z)$ whereas
$\bar{s}_{a_1, a_2, ...}$ stands for $S_{a_1, a_2, ...}(-1-z)$.

The constants are also written in a compact and  readable way  $\ln_2=\ln 2 \simeq 0.693147 $ and 
$ \text{LiHalf}_4=\text{Li}_4 \left(\frac{1}{2}\right)\simeq  0.517479$. For example, in this notation Eq.~(\ref{S1S21}) is written as Eq.~(\ref{s1s21}).

All other bilinear reflection identities are obtained by a trivial change of argument $z \leftrightarrow -1-z$.

\subsubsection{Reflection identities originating from $B_1 \otimes \bar{B}_3$} \label{app:B1B3}

\begin{eqnarray}
s_{-1} \bar{s}_{-3}& = & \frac{7 \zeta _2^2}{5}-\frac{s_2 \zeta _2}{2}-\frac{1}{2} \bar{s}_2
   \zeta _2-\ln _2 s_{-3}+\ln _2 s_3-\frac{3 \ln _2 \zeta _3}{2}-\frac{3}{4} s_{-1} \zeta _3
    \nonumber \\
   &&
   -\ln
   _2 \bar{s}_{-3}-\frac{3}{4} \zeta _3 \bar{s}_{-1}-\ln _2
   \bar{s}_3-s_{-3,-1}-\bar{s}_{-1,-3}
   \end{eqnarray} 
   
 \begin{eqnarray}
 s_{-1} \bar{s}_3& = & -\frac{3 \zeta
   _2^2}{5}+\frac{1}{2} s_{-2} \zeta _2-\frac{1}{2} \bar{s}_{-2} \zeta _2-\ln _2 s_{-3}+\ln _2
   s_3
   -\frac{3 \ln _2 \zeta _3}{2}
    \nonumber \\
   &&-\frac{3}{4} s_{-1} \zeta _3-\ln _2 \bar{s}_{-3}-\frac{3}{4}
   \zeta _3 \bar{s}_{-1}-\ln _2 \bar{s}_3+s_{3,-1}-\bar{s}_{-1,3}
   \end{eqnarray}

    \begin{eqnarray}
    s_{-1}
   \bar{s}_{-2,-1}& = & -\frac{\ln _2^4}{6}-s_{-2} \ln _2^2+s_2 \ln _2^2+\zeta _2 \ln
   _2^2-\bar{s}_{-2} \ln _2^2+\bar{s}_2 \ln _2^2-s_{-3} \ln _2+s_3 \ln _2
     \nonumber \\
   &&
   -6 \zeta _3 \ln
   _2
   -\bar{s}_{-3} \ln _2
   +\bar{s}_3 \ln _2-2 s_{-2,-1} \ln _2-2 \bar{s}_{-2,-1} \ln _2+\frac{33
   \zeta _2^2}{20}-4 \text{LiHalf}_4
   \nonumber \\
   &&
   +\frac{1}{2} s_{-2} \zeta _2
      +\frac{1}{4} s_{-1} \zeta
   _3
   +\frac{1}{2} \zeta _2 \bar{s}_{-2}+\frac{1}{4} \zeta _3
   \bar{s}_{-1}+s_{3,-1}+\bar{s}_{3,-1}-s_{-2,-1,-1}\nonumber \\
   &&
   -\bar{s}_{-2,-1,-1}-\bar{s}_{-1,-2,-1}
    \end{eqnarray}

       \begin{eqnarray}
   s_{-1} \bar{s}_{-2,1}& = & -\frac{\ln _2^4}{12}-\frac{1}{2} s_{-2} \ln _2^2+\frac{1}{2} s_2
   \ln _2^2+\frac{1}{2} \zeta _2 \ln _2^2-\frac{1}{2} \bar{s}_{-2} \ln _2^2+\frac{1}{2} \bar{s}_2
   \ln _2^2-s_{-3} \ln _2+s_3 \ln _2
   \nonumber \\
   &&
   -3 \zeta _3 \ln _2-s_{-2,-1} \ln _2+s_{-2,1} \ln
   _2-\bar{s}_{-2,-1} \ln _2-\bar{s}_{-2,1} \ln _2+\frac{61 \zeta _2^2}{40}-2
   \text{LiHalf}_4
   \nonumber \\
   &&
   -\frac{s_2 \zeta _2}{2}-\frac{5}{8} s_{-1} \zeta _3-\frac{1}{2} \zeta _2
   \bar{s}_{-2}-\frac{5}{8} \zeta _3 \bar{s}_{-1}-\frac{1}{2} \zeta _2
   \bar{s}_2-s_{-3,-1}+\bar{s}_{3,1}+s_{-2,1,-1}\nonumber \\
   &&
   -\bar{s}_{-2,-1,1}-\bar{s}_{-1,-2,1}
   \end{eqnarray}

   \begin{eqnarray}
   s_{-1} \bar{s}_{2,-1}& = & -\frac{\ln _2^4}{6}-s_{-2} \ln _2^2+s_2 \ln _2^2+\zeta _2 \ln
   _2^2+\bar{s}_{-2} \ln _2^2-\bar{s}_2 \ln _2^2-s_{-3} \ln _2
    \nonumber \\
   &&
   +s_3 \ln _2+\bar{s}_{-3} \ln
   _2-\bar{s}_3 \ln _2+2 s_{2,-1} \ln _2-2 \bar{s}_{2,-1} \ln _2+\frac{5 \zeta _2^2}{4}
    \nonumber \\
   &&-4
   \text{LiHalf}_4-\frac{s_2 \zeta _2}{2}+\frac{1}{4} s_{-1} \zeta _3+\frac{1}{4} \zeta _3
   \bar{s}_{-1}+\frac{1}{2} \zeta _2
   \bar{s}_2-s_{-3,-1}
    \nonumber \\
   &&+\bar{s}_{-3,-1}+s_{2,-1,-1}-\bar{s}_{-1,2,-1}-\bar{s}_{2,-1,-1}
   \end{eqnarray}

     \begin{eqnarray}
     s_{-1} \bar{s}_{2,1}& = & \frac{\ln _2^4}{6}-\frac{1}{2} s_{-2} \ln _2^2+\frac{1}{2} s_2 \ln
   _2^2-\zeta _2 \ln _2^2+\frac{1}{2} \bar{s}_{-2} \ln _2^2-\frac{1}{2} \bar{s}_2 \ln _2^2
    \nonumber \\
   &&
   -s_{-3}
   \ln _2+s_3 \ln _2+\frac{9 \zeta _3 \ln _2}{4}+s_{2,-1} \ln _2-s_{2,1} \ln _2-\bar{s}_{2,-1} \ln
   _2-\bar{s}_{2,1} \ln _2-\frac{7 \zeta _2^2}{4}
    \nonumber \\
   &&
   +4 \text{LiHalf}_4+\frac{1}{2} s_{-2} \zeta
   _2
   -\frac{5}{8} s_{-1} \zeta _3-\frac{1}{2} \zeta _2 \bar{s}_{-2}-\frac{5}{8} \zeta _3
   \bar{s}_{-1}-\frac{1}{2} \zeta _2
   \bar{s}_2
    \nonumber \\
   &&
   +s_{3,-1}+\bar{s}_{-3,1}-s_{2,1,-1}-\bar{s}_{-1,2,1}-\bar{s}_{2,-1,1}
   \end{eqnarray}

    \begin{eqnarray}
    s_{-1}\bar{s}_{-1,1,-1}& = & -\frac{7 \ln _2^4}{24}-s_{-2} \ln _2^2+s_2 \ln _2^2+\frac{9}{4} \zeta _2
   \ln _2^2-\frac{1}{2} \bar{s}_{-2} \ln _2^2+\frac{1}{2} \bar{s}_2 \ln _2^2
    \nonumber \\
   &&
   -\frac{1}{2} s_{-1,-1}
   \ln _2^2+\frac{3}{2} s_{-1,1} \ln _2^2+\frac{1}{2} \bar{s}_{-1,-1} \ln _2^2-\frac{1}{2}
   \bar{s}_{-1,1} \ln _2^2-s_{-3} \ln _2+s_3 \ln _2-3 \zeta _3 \ln _2
    \nonumber \\
   &&
   +\frac{3}{2} \zeta _2
   \bar{s}_{-1} \ln _2-2 s_{-2,-1} \ln _2-s_{-1,-2} \ln _2+s_{-1,2} \ln _2+\bar{s}_{-1,-2} \ln
   _2
    \nonumber \\
   &&
   -\bar{s}_{-1,2} \ln _2+2 s_{-1,1,-1} \ln _2
   -2 \bar{s}_{-1,1,-1} \ln _2+\frac{13 \zeta
   _2^2}{40}-\text{LiHalf}_4+\frac{1}{2} s_{-2} \zeta _2
    \nonumber \\
   &&
   +\frac{1}{4} s_{-1} \zeta _3-\frac{1}{8}
   \zeta _3 \bar{s}_{-1}-\frac{1}{2} \zeta _2 s_{-1,1}+s_{3,-1}+\frac{1}{2} \zeta _2
   \bar{s}_{-1,1}-s_{-2,-1,-1}
    \nonumber \\
   &&
   -s_{-1,-2,-1}+\bar{s}_{-1,-2,-1}+\bar{s}_{2,1,-1}+s_{-1,1,-1,-1}-2
   \bar{s}_{-1,-1,1,-1}-\bar{s}_{-1,1,-1,-1}
   \end{eqnarray}

    \begin{eqnarray}
    s_{-1} \bar{s}_{-1,1,1}& = & -\frac{\ln
   _2^4}{6}-\frac{1}{3} s_{-1} \ln _2^3-\frac{1}{2} s_{-2} \ln _2^2+\frac{1}{2} s_2 \ln
   _2^2
   +\frac{3}{4} \zeta _2 \ln _2^2-\frac{1}{2} s_{-1,-1} \ln _2^2
    \nonumber \\
   &&
   +\frac{1}{2} s_{-1,1} \ln
   _2^2+\frac{1}{2} \bar{s}_{-1,-1} \ln _2^2-\frac{1}{2} \bar{s}_{-1,1} \ln _2^2-s_{-3} \ln _2+s_3
   \ln _2+s_{-1} \zeta _2 \ln _2
    \nonumber \\
   &&
   -\frac{7 \zeta _3 \ln _2}{8}-s_{-2,-1} \ln _2+s_{-2,1} \ln
   _2-s_{-1,-2} \ln _2+s_{-1,2} \ln _2+s_{-1,1,-1} \ln _2
    \nonumber \\
   &&
   -s_{-1,1,1} \ln _2-\bar{s}_{-1,1,-1} \ln
   _2-\bar{s}_{-1,1,1} \ln _2-\frac{27 \zeta _2^2}{40}+4 \text{LiHalf}_4-\frac{s_2 \zeta
   _2}{2}
    \nonumber \\
   &&-\frac{9}{8} s_{-1} \zeta _3-\zeta _3 \bar{s}_{-1}-s_{-3,-1}+\frac{1}{2} \zeta _2
   s_{-1,-1}-\frac{1}{2} \zeta _2 \bar{s}_{-1,-1}-\frac{1}{2} \zeta _2
   \bar{s}_{-1,1}
    \nonumber \\
   &&
   +s_{-2,1,-1}+s_{-1,2,-1}+\bar{s}_{-1,-2,1}+\bar{s}_{2,1,1}-s_{-1,1,1,-1}
    \nonumber \\
   &&
   -2
   \bar{s}_{-1,-1,1,1}-\bar{s}_{-1,1,-1,1}
   \end{eqnarray}

     \begin{eqnarray}
     s_{-1} \bar{s}_{1,-2}& = & -\frac{\ln_2^4}{6}+\frac{3}{2} \zeta _2 \ln _2^2-s_{-3} \ln _2+s_3 \ln _2+\frac{1}{2} s_{-1} \zeta_2 \ln
   _2-\frac{1}{2} s_1 \zeta _2 \ln _2-\frac{15 \zeta _3 \ln _2}{4}
    \nonumber \\
   &&
   +\frac{1}{2} \zeta _2
   \bar{s}_{-1} \ln _2-\frac{1}{2} \zeta _2 \bar{s}_1 \ln _2+s_{1,-2} \ln _2-s_{1,2} \ln
   _2-\bar{s}_{1,-2} \ln _2-\bar{s}_{1,2} \ln _2+2 \zeta _2^2
    \nonumber \\
   &&
   -4 \text{LiHalf}_4-\frac{s_2 \zeta
   _2}{2}-\frac{1}{8} s_{-1} \zeta _3+\frac{5 s_1 \zeta _3}{8}+\frac{1}{2} \zeta _2
   \bar{s}_{-2}-\frac{1}{8} \zeta _3 \bar{s}_{-1}+\frac{13}{8} \zeta _3 \bar{s}_1
    \nonumber \\
   &&
   -\frac{1}{2}
   \zeta _2 \bar{s}_2-s_{-3,-1}-\frac{1}{2} \zeta _2 s_{1,-1}+\bar{s}_{-2,-2}-\frac{1}{2} \zeta _2
   \bar{s}_{1,-1}+s_{1,-2,-1} \nonumber \\
   &&
   -\bar{s}_{-1,1,-2}-\bar{s}_{1,-1,-2}
   \end{eqnarray}

    \begin{eqnarray}s_{-1}
   \bar{s}_{1,2}& = & -\frac{\ln _2^4}{12}+\zeta _2 \ln _2^2-s_{-3} \ln _2+s_3 \ln _2+\frac{1}{2}
   s_{-1} \zeta _2 \ln _2-\frac{1}{2} s_1 \zeta _2 \ln _2
    \nonumber \\
   &&
   -\frac{15 \zeta _3 \ln _2}{4}+\frac{1}{2}
   \zeta _2 \bar{s}_{-1} \ln _2-\frac{1}{2} \zeta _2 \bar{s}_1 \ln _2
   +s_{1,-2} \ln _2-s_{1,2} \ln
   _2
    \nonumber \\
   &&
   -\bar{s}_{1,-2} \ln _2-\bar{s}_{1,2} \ln _2
   +\frac{7 \zeta _2^2}{8}-2
   \text{LiHalf}_4+\frac{1}{2} s_{-2} \zeta _2-s_{-1} \zeta _3
    \nonumber \\
   &&
   +\frac{s_1 \zeta _3}{4}-\zeta _3
   \bar{s}_{-1}-\zeta _3 \bar{s}_1-\frac{1}{2} \zeta _2
   s_{1,-1}+s_{3,-1}+\bar{s}_{-2,2}
    \nonumber \\
   &&
   -\frac{1}{2} \zeta _2
   \bar{s}_{1,-1}-s_{1,2,-1}-\bar{s}_{-1,1,2}-\bar{s}_{1,-1,2}
   \end{eqnarray}

   \begin{eqnarray}
   s_{-1}\bar{s}_{1,-1,-1}& = & \frac{\ln _2^4}{8}-\frac{1}{6} s_{-1} \ln _2^3+\frac{7}{6} s_1 \ln
   _2^3-\frac{1}{6} \bar{s}_{-1} \ln _2^3+\frac{1}{6} \bar{s}_1 \ln _2^3-s_{-2} \ln _2^2
    \nonumber \\
   &&
   +s_2 \ln
   _2^2-\frac{1}{2} \zeta _2 \ln _2^2-\bar{s}_{-2} \ln _2^2+\bar{s}_2 \ln _2^2+2 s_{1,-1} \ln
   _2^2-s_{-3} \ln _2
    \nonumber \\
   &&
   +s_3 \ln _2-s_1 \zeta _2 \ln _2-3 \zeta _3 \ln _2-2 s_{-2,-1} \ln _2+s_{1,-2}
   \ln _2-s_{1,2} \ln _2
    \nonumber \\
   &&
   -\bar{s}_{1,-2} \ln _2+\bar{s}_{1,2} \ln _2
   +2 s_{1,-1,-1} \ln _2-2
   \bar{s}_{1,-1,-1} \ln _2+\frac{6 \zeta _2^2}{5}
    \nonumber \\
   &&
   -3 \text{LiHalf}_4+\frac{1}{2} s_{-2} \zeta
   _2-\frac{1}{4} s_{-1} \zeta _3+\frac{s_1 \zeta _3}{2}-\frac{1}{4} \zeta _3
   \bar{s}_{-1}-\frac{3}{4} \zeta _3 \bar{s}_1
    \nonumber \\
   &&
   -\frac{1}{2} \zeta _2 s_{1,-1}+s_{3,-1}+\frac{1}{2}
   \zeta _2
   \bar{s}_{1,-1}-s_{-2,-1,-1}-s_{1,2,-1}
    \nonumber \\
   &&
   +\bar{s}_{-2,-1,-1}+\bar{s}_{1,2,-1}+s_{1,-1,-1,-1}-\bar{
   s}_{-1,1,-1,-1}-2 \bar{s}_{1,-1,-1,-1}
   \end{eqnarray}

     \begin{eqnarray}
     s_{-1} \bar{s}_{1,-1,1}& = & -\frac{\ln
   _2^4}{12}-\frac{1}{6} s_{-1} \ln _2^3+\frac{5}{6} s_1 \ln _2^3-\frac{1}{6} \bar{s}_{-1} \ln
   _2^3+\frac{1}{6} \bar{s}_1 \ln _2^3
    \nonumber \\
   &&
   -\frac{1}{2} s_{-2} \ln _2^2+\frac{1}{2} s_2 \ln
   _2^2-\frac{1}{2} \zeta _2 \ln _2^2-\frac{1}{2} \bar{s}_{-2} \ln _2^2+\frac{1}{2} \bar{s}_2 \ln
   _2^2
    \nonumber \\
   &&+s_{1,-1} \ln _2^2-s_{-3} \ln _2+s_3 \ln _2-s_1 \zeta _2 \ln _2-\frac{19 \zeta _3 \ln
   _2}{8}
    \nonumber \\
   &&
   -\frac{3}{2} \zeta _2 \bar{s}_1 \ln _2-s_{-2,-1} \ln _2+s_{-2,1} \ln _2+s_{1,-2} \ln
   _2-s_{1,2} \ln _2
    \nonumber \\
   &&+s_{1,-1,-1} \ln _2-s_{1,-1,1} \ln _2-\bar{s}_{1,-1,-1} \ln
   _2-\bar{s}_{1,-1,1} \ln _2+\frac{49 \zeta _2^2}{20}
    \nonumber \\
   &&
   -6 \text{LiHalf}_4-\frac{s_2 \zeta
   _2}{2}+\frac{1}{8} s_{-1} \zeta _3+\frac{3 s_1 \zeta _3}{4}+\frac{1}{2} \zeta _2
   \bar{s}_{-2}+\frac{1}{8} \zeta _3 \bar{s}_{-1}
    \nonumber \\
   &&
   +\frac{3}{2} \zeta _3 \bar{s}_1-\frac{1}{2} \zeta
   _2 \bar{s}_2-s_{-3,-1}-\frac{1}{2} \zeta _2 s_{1,-1}-\zeta _2
   \bar{s}_{1,-1}+s_{-2,1,-1}
    \nonumber \\
   &&+s_{1,-2,-1}+\bar{s}_{-2,-1,1}+\bar{s}_{1,2,1}-s_{1,-1,1,-1}
    \nonumber \\
   &&
   -\bar{s}_
   {-1,1,-1,1}-2 \bar{s}_{1,-1,-1,1
   }\end{eqnarray}

    \begin{eqnarray}
    s_{-1} \bar{s}_{1,1,-1}& = & -\frac{\ln
   _2^4}{3}-\frac{1}{6} s_{-1} \ln _2^3-\frac{1}{6} s_1 \ln _2^3-\frac{1}{6} \bar{s}_{-1} \ln
   _2^3+\frac{1}{6} \bar{s}_1 \ln _2^3-s_{-2} \ln _2^2
    \nonumber \\
   &&
   +s_2 \ln _2^2+\frac{11}{4} \zeta _2 \ln
   _2^2+\frac{1}{2} \bar{s}_{-2} \ln _2^2-\frac{1}{2} \bar{s}_2 \ln _2^2+\frac{1}{2} s_{1,-1} \ln
   _2^2
    \nonumber \\
   &&
   -\frac{3}{2} s_{1,1} \ln _2^2+\frac{1}{2} \bar{s}_{1,-1} \ln _2^2-\frac{1}{2} \bar{s}_{1,1}
   \ln _2^2-s_{-3} \ln _2+s_3 \ln _2+\frac{1}{2} s_{-1} \zeta _2 \ln _2
    \nonumber \\
   &&
   +\frac{1}{2} s_1 \zeta _2
   \ln _2-4 \zeta _3 \ln _2+\frac{1}{2} \zeta _2 \bar{s}_{-1} \ln _2+\zeta _2 \bar{s}_1 \ln
   _2+s_{1,-2} \ln _2-s_{1,2} \ln _2
    \nonumber \\
   &&
   +2 s_{2,-1} \ln _2+\bar{s}_{1,-2} \ln _2-\bar{s}_{1,2} \ln
   _2-2 s_{1,1,-1} \ln _2-2 \bar{s}_{1,1,-1} \ln _2+\frac{8 \zeta _2^2}{5}
    \nonumber \\
   &&-4
   \text{LiHalf}_4-\frac{s_2 \zeta _2}{2}-\frac{s_1 \zeta _3}{4}-\frac{1}{8} \zeta _3
   \bar{s}_1-s_{-3,-1}+\frac{1}{2} \zeta _2 s_{1,1}+\frac{1}{2} \zeta _2
   \bar{s}_{1,1}
    \nonumber \\
   &&
   +s_{1,-2,-1}+s_{2,-1,-1}+\bar{s}_{-2,1,-1}+\bar{s}_{1,-2,-1}-s_{1,1,-1,-1}
    \nonumber \\
   &&
   -\bar{s}
   _{-1,1,1,-1}-\bar{s}_{1,-1,1,-1}-\bar{s}_{1,1,-1,-1}\end{eqnarray}

     \begin{eqnarray}
     s_{-1} \bar{s}_{1,1,1}& = &
   -\frac{\ln _2^4}{24}-\frac{1}{6} s_{-1} \ln _2^3+\frac{1}{6} s_1 \ln _2^3-\frac{1}{6}
   \bar{s}_{-1} \ln _2^3+\frac{1}{6} \bar{s}_1 \ln _2^3
   -\frac{1}{2} s_{-2} \ln _2^2
    \nonumber \\
   &&
   +\frac{1}{2}
   s_2 \ln _2^2+\frac{1}{4} \zeta _2 \ln _2^2+\frac{1}{2} s_{1,-1} \ln _2^2-\frac{1}{2} s_{1,1}
   \ln _2^2+\frac{1}{2} \bar{s}_{1,-1} \ln _2^2-\frac{1}{2} \bar{s}_{1,1} \ln _2^2
    \nonumber \\
   &&
   -s_{-3} \ln
   _2+s_3 \ln _2+\frac{1}{2} s_{-1} \zeta _2 \ln _2-\frac{1}{2} s_1 \zeta _2 \ln _2-\frac{7 \zeta
   _3 \ln _2}{8}
   +\frac{1}{2} \zeta _2 \bar{s}_{-1} \ln _2
    \nonumber \\
   &&
   -\frac{1}{2} \zeta _2 \bar{s}_1 \ln
   _2+s_{1,-2} \ln _2-s_{1,2} \ln _2+s_{2,-1} \ln _2-s_{2,1} \ln _2-s_{1,1,-1} \ln _2
    \nonumber \\
   &&
   +s_{1,1,1}
   \ln _2-\bar{s}_{1,1,-1} \ln _2
   -\bar{s}_{1,1,1} \ln _2-\frac{11 \zeta
   _2^2}{40}+\text{LiHalf}_4+\frac{1}{2} s_{-2} \zeta _2
    \nonumber \\
   &&
   -\frac{7}{8} s_{-1} \zeta _3+\frac{s_1
   \zeta _3}{4}-\frac{7}{8} \zeta _3 \bar{s}_{-1}-\frac{1}{8} \zeta _3 \bar{s}_1-\frac{1}{2} \zeta
   _2 s_{1,-1}+s_{3,-1}-\frac{1}{2} \zeta _2 \bar{s}_{1,-1}
    \nonumber \\
   &&
   -\frac{1}{2} \zeta _2
   \bar{s}_{1,1}-s_{1,2,-1}-s_{2,1,-1}+\bar{s}_{-2,1,1}+\bar{s}_{1,-2,1}+s_{1,1,1,-1}
    \nonumber \\
   &&
   -\bar{s}_{-1,
   1,1,1}-\bar{s}_{1,-1,1,1}-\bar{s}_{1,1,-1,1}
   \end{eqnarray}

    \begin{eqnarray}
    s_{-1} \bar{s}_{-1,-2}& = & \zeta _2 \ln
   _2^2-s_{-3} \ln _2+s_3 \ln _2+s_{-1} \zeta _2 \ln _2+\frac{3 \zeta _3 \ln _2}{2}-s_{-1,-2} \ln
   _2
    \nonumber \\
   &&
   +s_{-1,2} \ln _2-\bar{s}_{-1,-2} \ln _2-\bar{s}_{-1,2} \ln _2-\frac{11 \zeta
   _2^2}{10}+\frac{1}{2} s_{-2} \zeta _2-\frac{3}{4} s_{-1} \zeta _3
    \nonumber \\
   &&
   -\frac{1}{2} \zeta _2
   \bar{s}_{-2}+\frac{3}{2} \zeta _3 \bar{s}_{-1}+\frac{1}{2} \zeta _2 \bar{s}_2+\frac{1}{2} \zeta
   _2 s_{-1,-1}
   +s_{3,-1} \nonumber \\
   &&
   -\frac{1}{2} \zeta _2 \bar{s}_{-1,-1}+\bar{s}_{2,-2}-s_{-1,-2,-1}-2
   \bar{s}_{-1,-1,-2}
   \end{eqnarray}

    \begin{eqnarray}
    s_{-1} \bar{s}_{-1,2}& = & \frac{\ln _2^4}{3}-\zeta _2 \ln
   _2^2-s_{-3} \ln _2+s_3 \ln _2+s_{-1} \zeta _2 \ln _2+\frac{3 \zeta _3 \ln _2}{2}
    \nonumber \\
   &&
   -s_{-1,-2} \ln
   _2+s_{-1,2} \ln _2-\bar{s}_{-1,-2} \ln _2-\bar{s}_{-1,2} \ln _2-\frac{8 \zeta _2^2}{5}
    \nonumber \\
   &&
   +8
   \text{LiHalf}_4-\frac{s_2 \zeta _2}{2}-\frac{5}{4} s_{-1} \zeta _3-2 \zeta _3
   \bar{s}_{-1}-s_{-3,-1}
    \nonumber \\
   &&
   +\frac{1}{2} \zeta _2 s_{-1,-1}-\frac{1}{2} \zeta _2
   \bar{s}_{-1,-1}+\bar{s}_{2,2}+s_{-1,2,-1}-2 \bar{s}_{-1,-1,2}
   \end{eqnarray}

   \begin{eqnarray}s_{-1}
   \bar{s}_{-1,-1,-1}& = & -\frac{\ln _2^4}{2}-\frac{4}{3} s_{-1} \ln _2^3-s_{-2} \ln _2^2+s_2 \ln
   _2^2+\bar{s}_{-2} \ln _2^2-\bar{s}_2 \ln _2^2
    \nonumber \\
   &&
   -2 s_{-1,-1} \ln _2^2-s_{-3} \ln _2+s_3 \ln
   _2+s_{-1} \zeta _2 \ln _2+2 \zeta _3 \ln _2
    \nonumber \\
   &&
   -s_{-1,-2} \ln _2+s_{-1,2} \ln _2+2 s_{2,-1} \ln
   _2-\bar{s}_{-1,-2} \ln _2
    \nonumber \\
   &&
   +\bar{s}_{-1,2} \ln _2-2 s_{-1,-1,-1} \ln _2-2 \bar{s}_{-1,-1,-1} \ln
   _2-\frac{21 \zeta _2^2}{20}
    \nonumber \\
   &&
   +4 \text{LiHalf}_4-\frac{s_2 \zeta _2}{2}-\frac{3}{4} s_{-1} \zeta
   _3-\zeta _3 \bar{s}_{-1}-s_{-3,-1}
    \nonumber \\
   &&
   +\frac{1}{2} \zeta _2 s_{-1,-1}+\frac{1}{2} \zeta _2
   \bar{s}_{-1,-1}+s_{-1,2,-1}+s_{2,-1,-1}
    \nonumber \\
   &&
   +\bar{s}_{-1,2,-1}+\bar{s}_{2,-1,-1}-s_{-1,-1,-1,-1}-3
   \bar{s}_{-1,-1,-1,-1}
   \end{eqnarray}

    \begin{eqnarray}
    s_{-1} \bar{s}_{-1,-1,1}& = & -\frac{13 \ln _2^4}{24}-s_{-1} \ln
   _2^3-\frac{1}{2} s_{-2} \ln _2^2+\frac{1}{2} s_2 \ln _2^2
   +\frac{1}{2} \bar{s}_{-2} \ln
   _2^2-\frac{1}{2} \bar{s}_2 \ln _2^2
    \nonumber \\
   &&-s_{-1,-1} \ln _2^2-s_{-3} \ln _2+s_3 \ln _2+s_{-1} \zeta _2
   \ln _2+\frac{17 \zeta _3 \ln _2}{8}-\frac{3}{2} \zeta _2 \bar{s}_{-1} \ln _2
    \nonumber \\
   &&
   -s_{-1,-2} \ln
   _2+s_{-1,2} \ln _2+s_{2,-1} \ln _2-s_{2,1} \ln _2-s_{-1,-1,-1} \ln _2
    \nonumber \\
   &&
   +s_{-1,-1,1} \ln
   _2-\bar{s}_{-1,-1,-1} \ln _2-\bar{s}_{-1,-1,1} \ln _2
   -\frac{47 \zeta
   _2^2}{40}+\text{LiHalf}_4
    \nonumber \\
   &&
   +\frac{1}{2} s_{-2} \zeta _2-\frac{5}{8} s_{-1} \zeta _3-\frac{1}{2}
   \zeta _2 \bar{s}_{-2}+\frac{13}{8} \zeta _3 \bar{s}_{-1}+\frac{1}{2} \zeta _2
   \bar{s}_2+\frac{1}{2} \zeta _2 s_{-1,-1}
    \nonumber \\
   &&
   +s_{3,-1}-\zeta _2
   \bar{s}_{-1,-1}
   -s_{-1,-2,-1}-s_{2,1,-1}+\bar{s}_{-1,2,1}
    \nonumber \\
   &&
   +\bar{s}_{2,-1,1}+s_{-1,-1,1,-1}-3
   \bar{s}_{-1,-1,-1,1}
   \end{eqnarray}

    \begin{eqnarray}
    s_1 \bar{s}_{-3}& = & \frac{\ln _2^4}{6}-\zeta _2 \ln
   _2^2-\frac{37 \zeta _2^2}{20}+4 \text{LiHalf}_4+s_{-2} \zeta _2-\frac{7}{4} s_{-1} \zeta
   _3-\frac{3 s_1 \zeta _3}{4}-\zeta _2 \bar{s}_{-2}
    \nonumber \\
   &&
   -\frac{7}{4} \zeta _3 \bar{s}_{-1}+\frac{3}{4}
   \zeta _3 \bar{s}_1-s_{-3,1}+\bar{s}_{1,-3}
   \end{eqnarray}

    \begin{eqnarray}\label{s1s3}
    s_1 \bar{s}_3& = & \frac{8 \zeta
   _2^2}{5}-s_2 \zeta _2-\bar{s}_2 \zeta _2+s_1 \zeta _3-\zeta _3
   \bar{s}_1+s_{3,1}+\bar{s}_{1,3}
   \end{eqnarray}

    \begin{eqnarray}
    s_1 \bar{s}_{-2,-1}& = & -\frac{\ln
   _2^4}{4}-\frac{1}{2} s_{-2} \ln _2^2+\frac{1}{2} s_2 \ln _2^2+\frac{9}{2} \zeta _2 \ln
   _2^2-\frac{1}{2} \bar{s}_{-2} \ln _2^2+\frac{1}{2} \bar{s}_2 \ln _2^2+\frac{3}{2} s_{-1} \zeta
   _2 \ln _2
    \nonumber \\
   &&
   +\frac{3}{2} s_1 \zeta _2 \ln _2-3 \zeta _3 \ln _2-\bar{s}_{-3} \ln _2+\frac{3}{2}
   \zeta _2 \bar{s}_{-1} \ln _2
   -\frac{3}{2} \zeta _2 \bar{s}_1 \ln _2+\bar{s}_3 \ln _2-s_{-2,-1}
   \ln _2
    \nonumber \\
   &&
   -s_{-2,1} \ln _2-\bar{s}_{-2,-1} \ln _2+\bar{s}_{-2,1} \ln _2
   +\frac{15 \zeta _2^2}{8}-6
   \text{LiHalf}_4-\frac{1}{2} s_{-2} \zeta _2-\frac{s_2 \zeta _2}{2}
    \nonumber \\
   &&
   -\frac{5 s_1 \zeta
   _3}{8}+\frac{5}{8} \zeta _3 \bar{s}_1-\frac{1}{2} \zeta _2
   \bar{s}_2+s_{3,1}-\bar{s}_{-3,-1}-s_{-2,-1,1}+\bar{s}_{-2,1,-1}+\bar{s}_{1,-2,-1}
   \end{eqnarray}

    \begin{eqnarray}
    s_1   \bar{s}_{-2,1}& = & \frac{\ln _2^4}{4}-\frac{3}{2} \zeta _2 \ln _2^2-\frac{11 \zeta _2^2}{8}+6
   \text{LiHalf}_4+s_{-2} \zeta _2-\frac{21}{8} s_{-1} \zeta _3-\frac{5 s_1 \zeta _3}{8}
    \nonumber \\
   &&
   +\zeta _2
   \bar{s}_{-2}-\frac{21}{8} \zeta _3 \bar{s}_{-1}+\frac{5}{8} \zeta _3
   \bar{s}_1-s_{-3,1}-\bar{s}_{-3,1}+s_{-2,1,1}+\bar{s}_{-2,1,1}+\bar{s}_{1,-2,1}
   \end{eqnarray}

    \begin{eqnarray}
    s_1 \bar{s}_{2,-1}& = & -\frac{\ln _2^4}{12}-\frac{1}{2} s_{-2} \ln _2^2+\frac{1}{2} s_2 \ln
   _2^2-\frac{5}{2} \zeta _2 \ln _2^2+\frac{1}{2} \bar{s}_{-2} \ln _2^2-\frac{1}{2} \bar{s}_2 \ln
   _2^2-\frac{3}{2} s_{-1} \zeta _2 \ln _2
    \nonumber \\
   &&
   -\frac{3}{2} s_1 \zeta _2 \ln _2-\frac{9 \zeta _3 \ln
   _2}{4}+\bar{s}_{-3} \ln _2-\frac{3}{2} \zeta _2 \bar{s}_{-1} \ln _2+\frac{3}{2} \zeta _2
   \bar{s}_1 \ln _2-\bar{s}_3 \ln _2+s_{2,-1} \ln _2
    \nonumber \\
   &&
   +s_{2,1} \ln _2-\bar{s}_{2,-1} \ln
   _2+\bar{s}_{2,1} \ln _2+\frac{21 \zeta _2^2}{40}-2 \text{LiHalf}_4+\frac{1}{2} s_{-2} \zeta
   _2+\frac{s_2 \zeta _2}{2}+\frac{7}{8} s_{-1} \zeta _3
    \nonumber \\
   &&
   +\frac{s_1 \zeta _3}{4}-\frac{1}{2} \zeta
   _2 \bar{s}_{-2}+\frac{7}{8} \zeta _3 \bar{s}_{-1}
   -\frac{1}{4} \zeta _3
   \bar{s}_1-s_{-3,1}-\bar{s}_{3,-1}+s_{2,-1,1}
    \nonumber \\
   &&
   +\bar{s}_{1,2,-1}+\bar{s}_{2,1,-1}
   \end{eqnarray}

     \begin{eqnarray}\label{s1s21}
     s_1
   \bar{s}_{2,1}& = & \frac{6 \zeta _2^2}{5}-s_2 \zeta _2+\bar{s}_2 \zeta _2+2 s_1 \zeta _3-2 \zeta
   _3 \bar{s}_1+s_{3,1}-\bar{s}_{3,1}-s_{2,1,1}+\bar{s}_{1,2,1}+\bar{s}_{2,1,1}
   \end{eqnarray}

     \begin{eqnarray}
     s_1
   \bar{s}_{-1,1,-1}& = & -\frac{\ln _2^4}{4}+\frac{1}{6} s_{-1} \ln _2^3-\frac{1}{6} s_1 \ln
   _2^3-\frac{1}{6} \bar{s}_{-1} \ln _2^3+\frac{1}{6} \bar{s}_1 \ln _2^3-\frac{1}{2} s_{-2} \ln
   _2^2+\frac{1}{2} s_2 \ln _2^2
    \nonumber \\
   &&
   +\frac{3}{2} \zeta _2 \ln _2^2-\frac{1}{2} \bar{s}_{-2} \ln
   _2^2+\frac{1}{2} \bar{s}_2 \ln _2^2+s_{-1,1} \ln _2^2+\frac{5}{2} s_{-1} \zeta _2 \ln
   _2+\frac{1}{2} s_1 \zeta _2 \ln _2
    \nonumber \\
   &&
   -\frac{17 \zeta _3 \ln _2}{8}-\zeta _2 \bar{s}_{-1} \ln
   _2-\frac{1}{2} \zeta _2 \bar{s}_1 \ln _2-s_{-2,-1} \ln _2-s_{-2,1} \ln _2+\bar{s}_{-1,-2} \ln
   _2
    \nonumber \\
   &&
   -\bar{s}_{-1,2} \ln _2+s_{-1,1,-1} \ln _2+s_{-1,1,1} \ln _2-\bar{s}_{-1,1,-1} \ln
   _2+\bar{s}_{-1,1,1} \ln _2+\frac{19 \zeta _2^2}{20}
    \nonumber \\
   &&
   -2 \text{LiHalf}_4-\frac{1}{2} s_{-2} \zeta
   _2-\frac{s_2 \zeta _2}{2}-\frac{3}{4} s_{-1} \zeta _3+\frac{s_1 \zeta _3}{8}+\frac{1}{4} \zeta
   _3 \bar{s}_{-1}-\frac{1}{8} \zeta _3 \bar{s}_1
    \nonumber \\
   &&
   +\frac{1}{2} \zeta _2 s_{-1,-1}+\frac{1}{2} \zeta
   _2 s_{-1,1}+s_{3,1}-\frac{1}{2} \zeta _2
   \bar{s}_{-1,-1}-s_{-2,-1,1}-s_{-1,-2,1}
    \nonumber \\
   &&
   -\bar{s}_{-2,1,-1}-\bar{s}_{-1,2,-1}+s_{-1,1,-1,1}+2
   \bar{s}_{-1,1,1,-1}+\bar{s}_{1,-1,1,-1}
   \end{eqnarray}

     \begin{eqnarray}
     s_1 \bar{s}_{-1,1,1}& = & -\frac{\ln
   _2^4}{8}-\frac{1}{6} s_{-1} \ln _2^3-\frac{1}{6} s_1 \ln _2^3-\frac{1}{6} \bar{s}_{-1} \ln
   _2^3+\frac{1}{6} \bar{s}_1 \ln _2^3+\frac{1}{2} s_{-1} \zeta _2 \ln _2
    \nonumber \\
   &&
   +\frac{1}{2} s_1 \zeta _2
   \ln _2+\frac{1}{2} \zeta _2 \bar{s}_{-1} \ln _2-\frac{1}{2} \zeta _2 \bar{s}_1 \ln _2-\frac{13
   \zeta _2^2}{20}+3 \text{LiHalf}_4+s_{-2} \zeta _2
    \nonumber \\
   &&
   -\frac{23}{8} s_{-1} \zeta _3-\frac{7 s_1
   \zeta _3}{8}+\frac{1}{8} \zeta _3 \bar{s}_{-1}+\frac{7}{8} \zeta _3 \bar{s}_1-s_{-3,1}-\zeta _2
   s_{-1,1}
   +\zeta _2
   \bar{s}_{-1,1}+s_{-2,1,1}
    \nonumber \\
   &&
   +s_{-1,2,1}-\bar{s}_{-2,1,1}-\bar{s}_{-1,2,1}-s_{-1,1,1,1}+2
   \bar{s}_{-1,1,1,1}+\bar{s}_{1,-1,1,1}
   \end{eqnarray}

    \begin{eqnarray}
    s_1 \bar{s}_{1,-2}& = & \frac{\ln
   _2^4}{6}+\frac{1}{2} \zeta _2 \ln _2^2+\frac{3}{2} s_{-1} \zeta _2 \ln _2-\frac{3}{2} s_1 \zeta
   _2 \ln _2+\frac{3}{2} \zeta _2 \bar{s}_{-1} \ln _2-\frac{3}{2} \zeta _2 \bar{s}_1 \ln
   _2
    \nonumber \\
   &&
   -\frac{21 \zeta _2^2}{10}+4 \text{LiHalf}_4+s_{-2} \zeta _2-\frac{7}{4} s_{-1} \zeta
   _3+\frac{s_1 \zeta _3}{2}+\frac{1}{2} \zeta _2 \bar{s}_{-2}-\frac{7}{4} \zeta _3
   \bar{s}_{-1}
    \nonumber \\
   &&
   +\frac{1}{4} \zeta _3 \bar{s}_1-\frac{1}{2} \zeta _2 \bar{s}_2-s_{-3,1}-\frac{3}{2}
   \zeta _2 s_{1,-1}-\frac{1}{2} \zeta _2 s_{1,1}-\frac{3}{2} \zeta _2 \bar{s}_{1,-1}+\frac{1}{2}
   \zeta _2 \bar{s}_{1,1}
    \nonumber \\
   &&
   -\bar{s}_{2,-2}+s_{1,-2,1}+2 \bar{s}_{1,1,-2}
   \end{eqnarray}

   \begin{eqnarray}\label{s1s12}
   s_1 \bar{s}_{1,2}& = & \frac{6 \zeta _2^2}{5}-s_2 \zeta _2+s_{1,1} \zeta _2-\bar{s}_{1,1} \zeta
   _2+s_1 \zeta _3+2 \zeta _3 \bar{s}_1+s_{3,1}-\bar{s}_{2,2}-s_{1,2,1}+2
   \bar{s}_{1,1,2}
   \end{eqnarray}

   \begin{eqnarray}
   s_1 \bar{s}_{1,-1,-1}& = & \frac{\ln _2^4}{3}+s_1 \ln _2^3-\frac{1}{2}
   s_{-2} \ln _2^2+\frac{1}{2} s_2 \ln _2^2+\frac{11}{4} \zeta _2 \ln _2^2-\bar{s}_{-2} \ln
   _2^2+\bar{s}_2 \ln _2^2
    \nonumber \\
   &&
   +\frac{3}{2} s_{1,-1} \ln _2^2+\frac{1}{2} s_{1,1} \ln _2^2+\frac{1}{2}
   \bar{s}_{1,-1} \ln _2^2-\frac{1}{2} \bar{s}_{1,1} \ln _2^2+\frac{1}{2} s_{-1} \zeta _2 \ln
   _2 \nonumber \\
   &&
   +\frac{3}{2} s_1 \zeta _2 \ln _2-\frac{5 \zeta _3 \ln _2}{8}+\frac{1}{2} \zeta _2
   \bar{s}_{-1} \ln _2-\frac{1}{2} \zeta _2 \bar{s}_1 \ln _2-s_{-2,-1} \ln _2
    \nonumber \\
   &&
   -s_{-2,1} \ln
   _2-\bar{s}_{1,-2} \ln _2+\bar{s}_{1,2} \ln _2+s_{1,-1,-1} \ln _2+s_{1,-1,1} \ln
   _2-\bar{s}_{1,-1,-1} \ln _2
    \nonumber \\
   &&
   +\bar{s}_{1,-1,1} \ln _2+\frac{19 \zeta _2^2}{40}-\frac{1}{2} s_{-2}
   \zeta _2-\frac{s_2 \zeta _2}{2}-\frac{5 s_1 \zeta _3}{8}+\zeta _3 \bar{s}_1+\frac{1}{2} \zeta
   _2 s_{1,-1}
    \nonumber \\
   &&
   +\frac{1}{2} \zeta _2 s_{1,1}+s_{3,1}-\frac{1}{2} \zeta _2
   \bar{s}_{1,1}-s_{-2,-1,1}-s_{1,2,1}-\bar{s}_{1,-2,-1}-\bar{s}_{2,-1,-1}
    \nonumber \\
   &&
   +s_{1,-1,-1,1}+\bar{s}_{
   1,-1,1,-1}+2 \bar{s}_{1,1,-1,-1}
   \end{eqnarray}

    \begin{eqnarray}
    s_1 \bar{s}_{1,-1,1}& = & \frac{11 \ln
   _2^4}{24}+\frac{2}{3} s_1 \ln _2^3+\frac{3}{4} \zeta _2 \ln _2^2-\frac{1}{2} \bar{s}_{-2} \ln
   _2^2+\frac{1}{2} \bar{s}_2 \ln _2^2+\frac{1}{2} s_{1,-1} \ln _2^2
    \nonumber \\
   &&
   +\frac{1}{2} s_{1,1} \ln
   _2^2+\frac{1}{2} \bar{s}_{1,-1} \ln _2^2-\frac{1}{2} \bar{s}_{1,1} \ln _2^2+\frac{1}{2} s_{-1}
   \zeta _2 \ln _2-\frac{3}{2} s_1 \zeta _2 \ln _2
    \nonumber \\
   &&
   +\frac{1}{2} \zeta _2 \bar{s}_{-1} \ln _2+\zeta
   _2 \bar{s}_1 \ln _2-\frac{87 \zeta _2^2}{40}+5 \text{LiHalf}_4+s_{-2} \zeta _2-\frac{13}{8}
   s_{-1} \zeta _3
    \nonumber \\
   &&
   -\frac{s_1 \zeta _3}{8}+\frac{1}{2} \zeta _2 \bar{s}_{-2}-\frac{13}{8} \zeta _3
   \bar{s}_{-1}
   -\frac{1}{2} \zeta _2 \bar{s}_2-s_{-3,1}-\frac{3}{2} \zeta _2 s_{1,-1}
    \nonumber \\
   &&
   -\frac{1}{2}
   \zeta _2 s_{1,1}+\frac{1}{2} \zeta _2 \bar{s}_{1,-1}+\frac{1}{2} \zeta _2
   \bar{s}_{1,1}+s_{-2,1,1}+s_{1,-2,1}-\bar{s}_{1,-2,1}
    \nonumber \\
   &&
   -\bar{s}_{2,-1,1}-s_{1,-1,1,1}+\bar{s}_{1,-
   1,1,1}+2 \bar{s}_{1,1,-1,1}
   \end{eqnarray}

   \begin{eqnarray}
   s_1 \bar{s}_{1,1,-1}& = & -\frac{\ln _2^4}{8}-\frac{1}{3}
   s_1 \ln _2^3-\frac{1}{2} s_{-2} \ln _2^2+\frac{1}{2} s_2 \ln _2^2-\frac{3}{2} \zeta _2 \ln
   _2^2+\frac{1}{2} \bar{s}_{-2} \ln _2^2
    \nonumber \\
   &&-\frac{1}{2} \bar{s}_2 \ln _2^2-s_{1,1} \ln _2^2-2 s_1
   \zeta _2 \ln _2-\frac{25 \zeta _3 \ln _2}{8}-\frac{3}{2} \zeta _2 \bar{s}_1 \ln _2
    \nonumber \\
   &&
   +s_{2,-1} \ln
   _2+s_{2,1} \ln _2+\bar{s}_{1,-2} \ln _2-\bar{s}_{1,2} \ln _2-s_{1,1,-1} \ln _2-s_{1,1,1} \ln
   _2
    \nonumber \\
   &&-\bar{s}_{1,1,-1} \ln _2+\bar{s}_{1,1,1} \ln _2
   +\frac{13 \zeta
   _2^2}{40}-\text{LiHalf}_4+\frac{1}{2} s_{-2} \zeta _2+\frac{s_2 \zeta _2}{2}
    \nonumber \\
   &&
   -\frac{1}{8} s_{-1}
   \zeta _3+\frac{3 s_1 \zeta _3}{4}
   -\frac{1}{8} \zeta _3 \bar{s}_{-1}+\frac{1}{4} \zeta _3
   \bar{s}_1-s_{-3,1}-\frac{1}{2} \zeta _2 s_{1,-1}
    \nonumber \\
   &&
   -\frac{1}{2} \zeta _2 s_{1,1}-\frac{1}{2} \zeta
   _2 \bar{s}_{1,-1}+s_{1,-2,1}+s_{2,-1,1}-\bar{s}_{1,2,-1}
    \nonumber \\
   &&
   -\bar{s}_{2,1,-1}-s_{1,1,-1,1}+3
   \bar{s}_{1,1,1,-1}
   \end{eqnarray}

     \begin{eqnarray}
     s_1 \bar{s}_{1,1,1}& = & \frac{8 \zeta _2^2}{5}-s_2 \zeta _2+s_{1,1}
   \zeta _2+\bar{s}_{1,1} \zeta _2+2 s_1 \zeta _3+\zeta _3
   \bar{s}_1+s_{3,1}-s_{1,2,1}
   \nonumber \\
   &&
   -s_{2,1,1}-\bar{s}_{1,2,1}-\bar{s}_{2,1,1}+s_{1,1,1,1}+3
   \bar{s}_{1,1,1,1}
   \end{eqnarray}

    \begin{eqnarray}
    s_1 \bar{s}_{-1,-2}& = & \frac{\ln _2^4}{6}-\frac{1}{2} \zeta _2 \ln
   _2^2+2 s_{-1} \zeta _2 \ln _2-s_1 \zeta _2 \ln _2-\zeta _2 \bar{s}_{-1} \ln _2+\zeta _2
   \bar{s}_1 \ln _2
   \nonumber \\
   &&
   +4 \text{LiHalf}_4-s_2 \zeta _2-\frac{5}{8} s_{-1} \zeta _3+\frac{13 s_1 \zeta
   _3}{8}-\frac{1}{2} \zeta _2 \bar{s}_{-2}+\frac{1}{8} \zeta _3 \bar{s}_{-1}
   \nonumber \\
   &&
   -\frac{13}{8} \zeta
   _3 \bar{s}_1+\frac{1}{2} \zeta _2 \bar{s}_2+\frac{3}{2} \zeta _2 s_{-1,-1}+\frac{1}{2} \zeta _2
   s_{-1,1}+s_{3,1}-\bar{s}_{-2,-2}
   \nonumber \\
   &&
   -\frac{3}{2} \zeta _2 \bar{s}_{-1,-1}+\frac{1}{2} \zeta _2
   \bar{s}_{-1,1}-s_{-1,-2,1}+\bar{s}_{-1,1,-2}+\bar{s}_{1,-1,-2}
   \end{eqnarray}

   \begin{eqnarray}
   s_1 \bar{s}_{-1,2}& = &
   \frac{\ln _2^4}{12}+\frac{3}{2} \zeta _2 \ln _2^2+\frac{1}{2} s_{-1} \zeta _2 \ln
   _2+\frac{1}{2} s_1 \zeta _2 \ln _2+\frac{1}{2} \zeta _2 \bar{s}_{-1} \ln _2-\frac{1}{2} \zeta
   _2 \bar{s}_1 \ln _2
   \nonumber \\
   &&
   -\frac{9 \zeta _2^2}{8}+2 \text{LiHalf}_4+s_{-2} \zeta _2-2 s_{-1} \zeta
   _3-s_1 \zeta _3+\zeta _3 \bar{s}_{-1}+\zeta _3 \bar{s}_1-s_{-3,1}
   \nonumber \\
   &&
   -\zeta _2
   s_{-1,1}-\bar{s}_{-2,2}-\zeta _2
   \bar{s}_{-1,1}+s_{-1,2,1}+\bar{s}_{-1,1,2}+\bar{s}_{1,-1,2}
   \end{eqnarray}

    \begin{eqnarray}
    s_1 \bar{s}_{-1,-1,-1}& = &
   -\frac{5 \ln _2^4}{8}-\frac{7}{6} s_{-1} \ln _2^3-\frac{1}{6} s_1 \ln _2^3-\frac{1}{6}
   \bar{s}_{-1} \ln _2^3+\frac{1}{6} \bar{s}_1 \ln _2^3-\frac{1}{2} s_{-2} \ln _2^2
   \nonumber \\
   &&
   +\frac{1}{2}
   s_2 \ln _2^2-\frac{5}{4} \zeta _2 \ln _2^2+\bar{s}_{-2} \ln _2^2-\bar{s}_2 \ln _2^2-\frac{3}{2}
   s_{-1,-1} \ln _2^2-\frac{1}{2} s_{-1,1} \ln _2^2
   \nonumber \\
   &&
   +\frac{1}{2} \bar{s}_{-1,-1} \ln
   _2^2-\frac{1}{2} \bar{s}_{-1,1} \ln _2^2-\frac{3}{2} s_{-1} \zeta _2 \ln _2-\frac{1}{2} s_1
   \zeta _2 \ln _2-\frac{\zeta _3 \ln _2}{8}
   \nonumber \\
   &&
   -\frac{1}{2} \zeta _2 \bar{s}_{-1} \ln _2+\frac{1}{2}
   \zeta _2 \bar{s}_1 \ln _2+s_{2,-1} \ln _2+s_{2,1} \ln _2-\bar{s}_{-1,-2} \ln _2+\bar{s}_{-1,2}
   \ln _2
   \nonumber \\
   &&
   -s_{-1,-1,-1} \ln _2-s_{-1,-1,1} \ln _2-\bar{s}_{-1,-1,-1} \ln _2+\bar{s}_{-1,-1,1} \ln
   _2+\frac{3 \zeta _2^2}{40}
   \nonumber \\
   &&
   -\text{LiHalf}_4+\frac{1}{2} s_{-2} \zeta _2+\frac{s_2 \zeta
   _2}{2}+\frac{3}{8} s_{-1} \zeta _3-\frac{s_1 \zeta _3}{4}+\frac{3}{4} \zeta _3
   \bar{s}_{-1}+\frac{1}{4} \zeta _3 \bar{s}_1
   \nonumber \\
   &&
   -s_{-3,1}-\frac{1}{2} \zeta _2 s_{-1,-1}-\frac{1}{2}
   \zeta _2 s_{-1,1}-\frac{1}{2} \zeta _2
   \bar{s}_{-1,1}+s_{-1,2,1}+s_{2,-1,1}
   \nonumber \\
   &&
   -\bar{s}_{-2,-1,-1}-\bar{s}_{-1,-2,-1}-s_{-1,-1,-1,1}+\bar{
   s}_{-1,-1,1,-1}
   \nonumber \\
   &&
   +\bar{s}_{-1,1,-1,-1}+\bar{s}_{1,-1,-1,-1}
   \end{eqnarray}

     \begin{eqnarray}
     s_1 \bar{s}_{-1,-1,1}& = &
   -\frac{\ln _2^4}{2}-\frac{5}{6} s_{-1} \ln _2^3-\frac{1}{6} s_1 \ln _2^3-\frac{1}{6}
   \bar{s}_{-1} \ln _2^3+\frac{1}{6} \bar{s}_1 \ln _2^3+\frac{3}{4} \zeta _2 \ln _2^2
   \nonumber \\
   &&
   +\frac{1}{2}
   \bar{s}_{-2} \ln _2^2-\frac{1}{2} \bar{s}_2 \ln _2^2-\frac{1}{2} s_{-1,-1} \ln _2^2-\frac{1}{2}
   s_{-1,1} \ln _2^2+\frac{1}{2} \bar{s}_{-1,-1} \ln _2^2
   \nonumber \\
   &&
   -\frac{1}{2} \bar{s}_{-1,1} \ln
   _2^2+\frac{3}{2} s_{-1} \zeta _2 \ln _2-\frac{1}{2} s_1 \zeta _2 \ln _2+\zeta _2 \bar{s}_{-1}
   \ln _2+\frac{1}{2} \zeta _2 \bar{s}_1 \ln _2
   \nonumber \\
   &&
   +\frac{3 \zeta _2^2}{5}-s_2 \zeta _2+\frac{1}{4}
   s_{-1} \zeta _3+\frac{7 s_1 \zeta _3}{4}-\frac{1}{2} \zeta _2 \bar{s}_{-2}+\frac{1}{8} \zeta _3
   \bar{s}_{-1}-\frac{7}{4} \zeta _3 \bar{s}_1
   \nonumber \\
   &&
   +\frac{1}{2} \zeta _2 \bar{s}_2+\frac{3}{2} \zeta _2
   s_{-1,-1}+\frac{1}{2} \zeta _2 s_{-1,1}+s_{3,1}+\frac{1}{2} \zeta _2
   \bar{s}_{-1,-1}+\frac{1}{2} \zeta _2
   \bar{s}_{-1,1}
   \nonumber \\
   &&
   -s_{-1,-2,1}-s_{2,1,1}-\bar{s}_{-2,-1,1}-\bar{s}_{-1,-2,1}+s_{-1,-1,1,1}+\bar{s}_
   {-1,-1,1,1}
   \nonumber \\
   &&
   +\bar{s}_{-1,1,-1,1}+\bar{s}_{1,-1,-1,1} 
   \end{eqnarray}


\newpage

\subsubsection{Reflection identities originating from $B_2 \otimes \bar{B}_2$}

\begin{eqnarray}
 s_{-2} \bar{s}_{-2} & = &  \frac{8 \zeta _2^2}{5}-s_2 \zeta _2-\bar{s}_2 \zeta
   _2+s_{-2,-2}+\bar{s}_{-2,-2}
   \end{eqnarray}

    \begin{eqnarray}
    s_{-2} \bar{s}_2 & = &  \frac{\ln _2^4}{3}-2 \zeta _2 \ln
   _2^2-\frac{29 \zeta _2^2}{10}+8 \text{LiHalf}_4+\frac{1}{2} s_{-2} \zeta _2-\frac{s_2 \zeta
   _2}{2}-\frac{7}{2} s_{-1} \zeta _3-\frac{1}{2} \zeta _2 \bar{s}_{-2}
    \nonumber \\
   &&
   -\frac{7}{2} \zeta _3
   \bar{s}_{-1}-\frac{1}{2} \zeta _2 \bar{s}_2-s_{2,-2}+\bar{s}_{-2,2}
   \end{eqnarray}

    \begin{eqnarray}
    s_{-2}
   \bar{s}_{-1,1} & = &  \frac{\ln _2^4}{12}+\frac{1}{2} s_{-2} \ln _2^2+\frac{1}{2} s_2 \ln
   _2^2-\frac{1}{2} \bar{s}_{-2} \ln _2^2+\frac{1}{2} \bar{s}_2 \ln _2^2+\frac{3}{2} s_{-1} \zeta
   _2 \ln _2-\frac{3}{2} \zeta _2 \bar{s}_{-1} \ln _2
    \nonumber \\
   &&
   -\frac{11 \zeta _2^2}{40}+2
   \text{LiHalf}_4-s_2 \zeta _2-\frac{1}{8} s_{-1} \zeta _3+\frac{1}{2} \zeta _2
   \bar{s}_{-2}+\frac{5}{8} \zeta _3 \bar{s}_{-1}+\frac{1}{2} \zeta _2
   \bar{s}_2
    \nonumber \\
   &&
   +s_{-2,-2}+\frac{3}{2} \zeta _2 s_{-1,-1}-\frac{1}{2} \zeta _2 s_{-1,1}-\frac{3}{2}
   \zeta _2 \bar{s}_{-1,-1}-\frac{1}{2} \zeta _2
   \bar{s}_{-1,1}-\bar{s}_{3,1}
    \nonumber \\
   &&
   -s_{-1,1,-2}+\bar{s}_{-2,-1,1}+\bar{s}_{-1,-2,1}
   \end{eqnarray}

    \begin{eqnarray}
       s_{-2}
   \bar{s}_{1,-1} & = &  -\frac{\ln _2^4}{4}-\frac{1}{2} s_{-2} \ln _2^2-\frac{1}{2} s_2 \ln _2^2+4
   \zeta _2 \ln _2^2+\frac{1}{2} \bar{s}_{-2} \ln _2^2-\frac{1}{2} \bar{s}_2 \ln _2^2
    \nonumber \\
   &&
   +\frac{3}{2}
   s_{-1} \zeta _2 \ln _2-s_1 \zeta _2 \ln _2-\frac{15 \zeta _3 \ln _2}{4}-\bar{s}_{-3} \ln
   _2+\frac{3}{2} \zeta _2 \bar{s}_{-1} \ln _2
    \nonumber \\
   &&
   +\bar{s}_3 \ln _2-s_{1,-2} \ln _2-s_{1,2} \ln
   _2+\bar{s}_{1,-2} \ln _2-\bar{s}_{1,2} \ln _2+\frac{17 \zeta _2^2}{8}
    \nonumber \\
   &&
   -6
   \text{LiHalf}_4+\frac{1}{2} s_{-2} \zeta _2-\frac{s_2 \zeta _2}{2}+\frac{13 s_1 \zeta
   _3}{8}+\frac{5}{8} \zeta _3 \bar{s}_1-\frac{1}{2} \zeta _2 \bar{s}_2
    \nonumber \\
   &&
   +s_{-2,-2}-\frac{1}{2}
   \zeta _2 s_{1,-1}-\bar{s}_{-3,-1}-\frac{1}{2} \zeta _2
   \bar{s}_{1,-1}-s_{1,-1,-2}
    \nonumber \\
   &&
   +\bar{s}_{-2,1,-1}+\bar{s}_{1,-2,-1}
   \end{eqnarray}

     \begin{eqnarray}
     s_{-2}
   \bar{s}_{1,1} & = &  \frac{\ln _2^4}{4}+\frac{3}{2} s_{-1} \zeta _2 \ln _2-\frac{3}{2} s_1 \zeta _2
   \ln _2+\frac{3}{2} \zeta _2 \bar{s}_{-1} \ln _2-\frac{3}{2} \zeta _2 \bar{s}_1 \ln _2-\frac{13
   \zeta _2^2}{8}
    \nonumber \\
   &&
   +6 \text{LiHalf}_4+\frac{1}{2} s_{-2} \zeta _2-\frac{s_2 \zeta
   _2}{2}-\frac{21}{8} s_{-1} \zeta _3+\frac{s_1 \zeta _3}{8}+\zeta _2 \bar{s}_{-2}
    \nonumber \\
   &&
   -\frac{21}{8}
   \zeta _3 \bar{s}_{-1}+\frac{5}{8} \zeta _3 \bar{s}_1-\frac{3}{2} \zeta _2 s_{1,-1}+\frac{1}{2}
   \zeta _2 s_{1,1}-s_{2,-2}-\bar{s}_{-3,1}
    \nonumber \\
   &&
   -\frac{3}{2} \zeta _2 \bar{s}_{1,-1}-\frac{1}{2} \zeta
   _2 \bar{s}_{1,1}+s_{1,1,-2}+\bar{s}_{-2,1,1}+\bar{s}_{1,-2,1}
   \end{eqnarray}

    \begin{eqnarray}
    s_{-2}
   \bar{s}_{-1,-1} & = &  \frac{\ln _2^4}{6}+s_{-2} \ln _2^2+s_2 \ln _2^2-\bar{s}_{-2} \ln
   _2^2+\bar{s}_2 \ln _2^2+s_{-1} \zeta _2 \ln _2-\frac{3 \zeta _3 \ln _2}{2}
    \nonumber \\
   &&
   +\bar{s}_{-3} \ln
   _2-\bar{s}_3 \ln _2+s_{-1,-2} \ln _2+s_{-1,2} \ln _2+\bar{s}_{-1,-2} \ln _2-\bar{s}_{-1,2} \ln
   _2
    \nonumber \\
   &&
   -\frac{23 \zeta _2^2}{20}+4 \text{LiHalf}_4+\frac{1}{2} s_{-2} \zeta _2-\frac{s_2 \zeta
   _2}{2}-\frac{5}{2} s_{-1} \zeta _3-\frac{1}{2} \zeta _2 \bar{s}_{-2}
    \nonumber \\
   &&
   -\frac{1}{4} \zeta _3
   \bar{s}_{-1}+\frac{1}{2} \zeta _2 s_{-1,-1}-s_{2,-2}-\frac{1}{2} \zeta _2
   \bar{s}_{-1,-1}-\bar{s}_{3,-1}+s_{-1,-1,-2}
    \nonumber \\
   &&
   +\bar{s}_{-2,-1,-1}+\bar{s}_{-1,-2,-1}
   \end{eqnarray}

    \begin{eqnarray}
    s_2
   \bar{s}_2 & = &  \frac{12 \zeta _2^2}{5}-s_{2,2}-\bar{s}_{2,2}
   \end{eqnarray}

   \begin{eqnarray}
   s_2 \bar{s}_{-1,1} & = & 
   \frac{\ln _2^4}{6}+\frac{1}{2} s_{-2} \ln _2^2+\frac{1}{2} s_2 \ln _2^2+\zeta _2 \ln
   _2^2-\frac{1}{2} \bar{s}_{-2} \ln _2^2+\frac{1}{2} \bar{s}_2 \ln _2^2+\frac{3}{2} s_{-1} \zeta
   _2 \ln _2
    \nonumber \\
   &&
   +\frac{3}{2} \zeta _2 \bar{s}_{-1} \ln _2-\frac{5 \zeta _2^2}{4}+4
   \text{LiHalf}_4-\frac{1}{2} s_{-2} \zeta _2-\frac{s_2 \zeta _2}{2}-\frac{29}{8} s_{-1} \zeta
   _3-\frac{1}{2} \zeta _2 \bar{s}_{-2}
    \nonumber \\
   &&
   -\frac{5}{8} \zeta _3 \bar{s}_{-1}-\frac{1}{2} \zeta _2
   \bar{s}_2+s_{-2,2}+\zeta _2 s_{-1,1}+\bar{s}_{-3,1}+\zeta _2
   \bar{s}_{-1,1}-s_{-1,1,2}
    \nonumber \\
   &&
   -\bar{s}_{-1,2,1}-\bar{s}_{2,-1,1}
   \end{eqnarray}

    \begin{eqnarray}
    s_2 \bar{s}_{1,-1} & = & 
   \frac{\ln _2^4}{12}-\frac{1}{2} s_{-2} \ln _2^2-\frac{1}{2} s_2 \ln _2^2-\zeta _2 \ln
   _2^2+\frac{1}{2} \bar{s}_{-2} \ln _2^2-\frac{1}{2} \bar{s}_2 \ln _2^2+\frac{1}{2} s_1 \zeta _2
   \ln _2-\frac{15 \zeta _3 \ln _2}{4}
    \nonumber \\
   &&
   -\bar{s}_{-3} \ln _2-\frac{3}{2} \zeta _2 \bar{s}_1 \ln
   _2+\bar{s}_3 \ln _2-s_{1,-2} \ln _2-s_{1,2} \ln _2+\bar{s}_{1,-2} \ln _2-\bar{s}_{1,2} \ln
   _2-\frac{11 \zeta _2^2}{40}
    \nonumber \\
   &&
   +2 \text{LiHalf}_4-\frac{7}{8} s_{-1} \zeta _3-s_1 \zeta
   _3+\frac{1}{2} \zeta _2 \bar{s}_{-2}-\frac{7}{8} \zeta _3 \bar{s}_{-1}+\frac{1}{4} \zeta _3
   \bar{s}_1+s_{-2,2}-\frac{1}{2} \zeta _2 s_{1,-1}
    \nonumber \\
   &&
   -\frac{1}{2} \zeta _2
   \bar{s}_{1,-1}+\bar{s}_{3,-1}-s_{1,-1,2}-\bar{s}_{1,2,-1}-\bar{s}_{2,1,-1}
   \end{eqnarray}

    \begin{eqnarray}
    s_2
   \bar{s}_{1,1} & = &  \frac{4 \zeta _2^2}{5}-\bar{s}_2 \zeta _2-s_{1,1} \zeta _2+\bar{s}_{1,1} \zeta
   _2+s_1 \zeta _3+2 \zeta _3
   \bar{s}_1-s_{2,2}+\bar{s}_{3,1}+s_{1,1,2}
    \nonumber \\
   &&
   -\bar{s}_{1,2,1}-\bar{s}_{2,1,1}
   \end{eqnarray}

    \begin{eqnarray}
    s_2
   \bar{s}_{-1,-1} & = &  -\frac{\ln _2^4}{6}+s_{-2} \ln _2^2+s_2 \ln _2^2+2 \zeta _2 \ln
   _2^2-\bar{s}_{-2} \ln _2^2+\bar{s}_2 \ln _2^2+s_{-1} \zeta _2 \ln _2-\frac{3 \zeta _3 \ln
   _2}{2}
    \nonumber \\
   &&
   +\bar{s}_{-3} \ln _2-\bar{s}_3 \ln _2+s_{-1,-2} \ln _2+s_{-1,2} \ln _2+\bar{s}_{-1,-2}
   \ln _2-\bar{s}_{-1,2} \ln _2+\frac{7 \zeta _2^2}{4}
    \nonumber \\
   &&
   -4 \text{LiHalf}_4+s_{-1} \zeta
   _3+\frac{1}{4} \zeta _3 \bar{s}_{-1}+\frac{1}{2} \zeta _2 \bar{s}_2+\frac{1}{2} \zeta _2
   s_{-1,-1}-s_{2,2}+\bar{s}_{-3,-1}
    \nonumber \\
   &&
   -\frac{1}{2} \zeta _2
   \bar{s}_{-1,-1}+s_{-1,-1,2}-\bar{s}_{-1,2,-1}-\bar{s}_{2,-1,-1}
   \end{eqnarray}

   \begin{eqnarray}
   s_{-1,1}
   \bar{s}_{-1,1} & = &  -\frac{\ln _2^4}{4}-\frac{1}{2} s_{-2} \ln _2^2+\frac{1}{2} s_2 \ln
   _2^2-\frac{1}{2} \zeta _2 \ln _2^2-\frac{1}{2} \bar{s}_{-2} \ln _2^2+\frac{1}{2} \bar{s}_2 \ln
   _2^2+s_{-1,1} \ln _2^2
    \nonumber \\
   &&
   +\bar{s}_{-1,1} \ln _2^2-\frac{3}{2} s_{-1} \zeta _2 \ln _2-\frac{3}{2}
   \zeta _2 \bar{s}_{-1} \ln _2-\frac{19 \zeta _2^2}{20}+6 \text{LiHalf}_4+\frac{1}{2} s_{-2}
   \zeta _2+\frac{s_2 \zeta _2}{2}
    \nonumber \\
   &&
   -\frac{3}{8} s_{-1} \zeta _3+\frac{1}{2} \zeta _2
   \bar{s}_{-2}-\frac{3}{8} \zeta _3 \bar{s}_{-1}+\frac{1}{2} \zeta _2 \bar{s}_2-2 \zeta _2
   s_{-1,-1}-\zeta _2 s_{-1,1}-s_{3,1}
    \nonumber \\
   &&
   -2 \zeta _2 \bar{s}_{-1,-1}-\zeta _2
   \bar{s}_{-1,1}-\bar{s}_{3,1}+s_{-2,-1,1}+2 s_{-1,-2,1}+s_{2,1,1}+\bar{s}_{-2,-1,1}
    \nonumber \\
   &&
   +2
   \bar{s}_{-1,-2,1}+\bar{s}_{2,1,1}-2 s_{-1,-1,1,1}-s_{-1,1,-1,1}-2
   \bar{s}_{-1,-1,1,1}-\bar{s}_{-1,1,-1,1}
   \end{eqnarray}

    \begin{eqnarray}
    s_{-1,1} \bar{s}_{1,-1} & = &  -\frac{\ln
   _2^4}{3}-\frac{2}{3} s_{-1} \ln _2^3+\frac{2}{3} s_1 \ln _2^3-\frac{1}{3} \bar{s}_{-1} \ln
   _2^3+\frac{1}{3} \bar{s}_1 \ln _2^3-\frac{1}{2} \zeta _2 \ln _2^2-s_{-1,1} \ln
   _2^2
    \nonumber \\
   &&
   +\bar{s}_{1,-1} \ln _2^2-\frac{3}{2} s_{-1} \zeta _2 \ln _2-2 s_1 \zeta _2 \ln
   _2-\frac{\zeta _3 \ln _2}{4}-\bar{s}_{-3} \ln _2+2 \zeta _2 \bar{s}_{-1} \ln _2
    \nonumber \\
   &&
   -\frac{1}{2}
   \zeta _2 \bar{s}_1 \ln _2+\bar{s}_3 \ln _2+s_{-2,-1} \ln _2+s_{-2,1} \ln _2-\bar{s}_{-2,-1} \ln
   _2+\bar{s}_{-2,1} \ln _2
    \nonumber \\
   &&
   -\bar{s}_{-1,-2} \ln _2+\bar{s}_{-1,2} \ln _2+\bar{s}_{1,-2} \ln
   _2-\bar{s}_{1,2} \ln _2-s_{-1,1,-1} \ln _2-s_{-1,1,1} \ln _2
    \nonumber \\
   &&
   -s_{1,-1,-1} \ln _2-s_{1,-1,1} \ln
   _2+\bar{s}_{-1,1,-1} \ln _2-\bar{s}_{-1,1,1} \ln _2+\bar{s}_{1,-1,-1} \ln _2
   -\bar{s}_{1,-1,1}
   \ln _2
    \nonumber \\
   &&
   +\frac{47 \zeta _2^2}{40}-4 \text{LiHalf}_4+s_{-2} \zeta _2+\frac{3}{4} s_{-1} \zeta
   _3+\frac{3 s_1 \zeta _3}{2}-\frac{1}{4} \zeta _3 \bar{s}_{-1}+\frac{3}{4} \zeta _3
   \bar{s}_1
    \nonumber \\
   &&
   -\frac{1}{2} \zeta _2 \bar{s}_2-\frac{1}{2} \zeta _2 s_{-1,-1}-\frac{1}{2} \zeta _2
   s_{-1,1}-\zeta _2 s_{1,-1}-s_{3,1}-\bar{s}_{-3,-1}+\frac{1}{2} \zeta _2
   \bar{s}_{-1,-1}
    \nonumber \\
   &&
   -\frac{1}{2} \zeta _2 \bar{s}_{1,-1}+2 s_{-2,-1,1}+s_{-1,-2,1}+s_{1,2,1}+2
   \bar{s}_{-2,1,-1}+\bar{s}_{-1,2,-1}
    \nonumber \\
   &&
   +\bar{s}_{1,-2,-1}-s_{-1,1,-1,1}-2 s_{1,-1,-1,1}-2
   \bar{s}_{-1,1,1,-1}-\bar{s}_{1,-1,1,-1}
   \end{eqnarray}

   \begin{eqnarray}
   s_{-1,1} \bar{s}_{1,1} & = &  -\frac{\ln
   _2^4}{8}-\frac{1}{3} s_{-1} \ln _2^3+\frac{1}{3} s_1 \ln _2^3-\frac{1}{3} \bar{s}_{-1} \ln
   _2^3+\frac{1}{3} \bar{s}_1 \ln _2^3-\frac{1}{2} s_{-2} \ln _2^2
   +\frac{1}{2} s_2 \ln
   _2^2
    \nonumber \\
   &&
   +\frac{3}{4} \zeta _2 \ln _2^2
   +\frac{1}{2} s_{1,-1} \ln _2^2-\frac{1}{2} s_{1,1} \ln
   _2^2+\frac{1}{2} \bar{s}_{1,-1} \ln _2^2+\frac{1}{2} \bar{s}_{1,1} \ln _2^2+s_{-1} \zeta _2 \ln
   _2
    \nonumber \\
   &&
   +\frac{1}{2} s_1 \zeta _2 \ln _2+\zeta _2 \bar{s}_{-1} \ln _2-\zeta _2 \bar{s}_1 \ln
   _2-\frac{39 \zeta _2^2}{40}+3 \text{LiHalf}_4-\frac{1}{2} s_{-2} \zeta _2-\frac{s_2 \zeta
   _2}{2}
    \nonumber \\
   &&
   +\frac{1}{4} s_{-1} \zeta _3+\frac{s_1 \zeta _3}{8}+\zeta _2 \bar{s}_{-2}-\frac{11}{4}
   \zeta _3 \bar{s}_{-1}-\frac{1}{4} \zeta _3 \bar{s}_1+s_{-3,1}+\zeta _2 s_{-1,1}+\frac{1}{2}
   \zeta _2 s_{1,-1}
    \nonumber \\
   &&
   +\frac{1}{2} \zeta _2 s_{1,1}-\bar{s}_{-3,1}-\zeta _2
   \bar{s}_{-1,1}
   -\frac{3}{2} \zeta _2 \bar{s}_{1,-1}-\frac{1}{2} \zeta _2
   \bar{s}_{1,1}-s_{-2,1,1}-s_{-1,2,1}
    \nonumber \\
   &&
   -s_{1,-2,1}-s_{2,-1,1}+2
   \bar{s}_{-2,1,1}+\bar{s}_{-1,2,1}+\bar{s}_{1,-2,1}+s_{-1,1,1,1}+s_{1,-1,1,1}
    \nonumber \\
   &&
   +s_{1,1,-1,1}-2
   \bar{s}_{-1,1,1,1}-\bar{s}_{1,-1,1,1}
   \end{eqnarray}

    \begin{eqnarray}
    s_{-1,1} \bar{s}_{-1,-1} & = &  -\frac{11 \ln
   _2^4}{24}+\frac{11}{4} \zeta _2 \ln _2^2-\frac{1}{2} \bar{s}_{-2} \ln _2^2+\frac{1}{2}
   \bar{s}_2 \ln _2^2+\frac{3}{2} s_{-1,-1} \ln _2^2
   +\frac{1}{2} s_{-1,1} \ln _2^2
    \nonumber \\
   &&
   +\frac{1}{2}
   \bar{s}_{-1,-1} \ln _2^2+\frac{1}{2} \bar{s}_{-1,1} \ln _2^2+\frac{7}{2} s_{-1} \zeta _2 \ln
   _2-2 \zeta _3 \ln _2+\bar{s}_{-3} \ln _2-\bar{s}_3 \ln _2
    \nonumber \\
   &&-s_{2,-1} \ln _2-s_{2,1} \ln _2+2
   \bar{s}_{-1,-2} \ln _2-2 \bar{s}_{-1,2} \ln _2-\bar{s}_{2,-1} \ln _2+\bar{s}_{2,1} \ln _2
    \nonumber \\
   &&
   +2
   s_{-1,-1,-1} \ln _2+2 s_{-1,-1,1} \ln _2+2 \bar{s}_{-1,-1,-1} \ln _2-2 \bar{s}_{-1,-1,1} \ln
   _2-\frac{33 \zeta _2^2}{40}
    \nonumber \\
   &&
   +3 \text{LiHalf}_4-s_2 \zeta _2-\frac{9}{4} s_{-1} \zeta
   _3-\frac{1}{2} \zeta _2 \bar{s}_{-2}-\frac{3}{8} \zeta _3 \bar{s}_{-1}+s_{-3,1}+\frac{3}{2}
   \zeta _2 s_{-1,-1}
    \nonumber \\
   &&
   +\frac{1}{2} \zeta _2 s_{-1,1}-\frac{1}{2} \zeta _2
   \bar{s}_{-1,-1}+\frac{1}{2} \zeta _2 \bar{s}_{-1,1}
   -\bar{s}_{3,-1}-2 s_{-1,2,1}-2
   s_{2,-1,1}
    \nonumber \\
   &&
   +\bar{s}_{-2,-1,-1}+2 \bar{s}_{-1,-2,-1} 
   +\bar{s}_{2,1,-1}+3 s_{-1,-1,-1,1}
    \nonumber \\
   &&
   -2
   \bar{s}_{-1,-1,1,-1}-\bar{s}_{-1,1,-1,-1}
   \end{eqnarray}

   \begin{eqnarray}
   s_{1,-1} \bar{s}_{1,-1} & = &  -\frac{\ln
   _2^4}{12}-\frac{1}{2} s_{-2} \ln _2^2+\frac{1}{2} s_2 \ln _2^2+\frac{15}{2} \zeta _2 \ln
   _2^2-\frac{1}{2} \bar{s}_{-2} \ln _2^2+\frac{1}{2} \bar{s}_2 \ln _2^2
    \nonumber \\
   &&
   +s_{1,-1} \ln _2^2-2
   s_{1,1} \ln _2^2+\bar{s}_{1,-1} \ln _2^2-2 \bar{s}_{1,1} \ln _2^2-s_{-3} \ln _2+s_3 \ln
   _2
    \nonumber \\
   &&
   +s_{-1} \zeta _2 \ln _2+\frac{3}{2} s_1 \zeta _2 \ln _2-8 \zeta _3 \ln _2-\bar{s}_{-3} \ln
   _2+\zeta _2 \bar{s}_{-1} \ln _2+\frac{3}{2} \zeta _2 \bar{s}_1 \ln _2
   \nonumber \\
   &&
   +\bar{s}_3 \ln _2+2
   s_{1,-2} \ln _2-2 s_{1,2} \ln _2+2 s_{2,-1} \ln _2+2 \bar{s}_{1,-2} \ln _2-2 \bar{s}_{1,2} \ln
   _2
   \nonumber \\
   &&
   +2 \bar{s}_{2,-1} \ln _2-4 s_{1,1,-1} \ln _2-4 \bar{s}_{1,1,-1} \ln _2+\frac{12 \zeta
   _2^2}{5}-6 \text{LiHalf}_4-\frac{s_2 \zeta _2}{2}
   \nonumber \\
   &&
   -\frac{3 s_1 \zeta _3}{8}-\frac{3}{8} \zeta _3
   \bar{s}_1-\frac{1}{2} \zeta _2 \bar{s}_2-s_{-3,-1}+\zeta _2 s_{1,1}-\bar{s}_{-3,-1}+\zeta _2
   \bar{s}_{1,1}
   \nonumber \\
   &&
   +s_{-2,1,-1}+2 s_{1,-2,-1}+s_{2,-1,-1}+\bar{s}_{-2,1,-1}+2
   \bar{s}_{1,-2,-1}+\bar{s}_{2,-1,-1}
   \nonumber \\
   &&
   -s_{1,-1,1,-1}-2 s_{1,1,-1,-1}-\bar{s}_{1,-1,1,-1}-2
   \bar{s}_{1,1,-1,-1}
   \end{eqnarray}

   \begin{eqnarray}
   s_{1,-1} \bar{s}_{1,1} & = &  \frac{\ln _2^4}{8}-\frac{1}{4} \zeta _2
   \ln _2^2-\frac{1}{2} \bar{s}_{-2} \ln _2^2+\frac{1}{2} \bar{s}_2 \ln _2^2
   +\frac{1}{2} s_{1,-1}
   \ln _2^2-\frac{1}{2} s_{1,1} \ln _2^2+\frac{1}{2} \bar{s}_{1,-1} \ln _2^2
   \nonumber \\
   &&
   -\frac{3}{2}
   \bar{s}_{1,1} \ln _2^2-s_{-3} \ln _2+s_3 \ln _2+\frac{1}{2} s_{-1} \zeta _2 \ln _2-2 s_1 \zeta
   _2 \ln _2-4 \zeta _3 \ln _2+\frac{1}{2} \zeta _2 \bar{s}_{-1} \ln _2
   \nonumber \\
   &&
   -\frac{5}{2} \zeta _2
   \bar{s}_1 \ln _2+2 s_{1,-2} \ln _2-2 s_{1,2} \ln _2+s_{2,-1} \ln _2-s_{2,1} \ln
   _2+\bar{s}_{2,-1} \ln _2+\bar{s}_{2,1} \ln _2
   \nonumber \\
   &&
   -2 s_{1,1,-1} \ln _2
   +2 s_{1,1,1} \ln _2-2
   \bar{s}_{1,1,-1} \ln _2
   -2 \bar{s}_{1,1,1} \ln _2+\frac{\zeta
   _2^2}{20}+\text{LiHalf}_4+\frac{1}{2} s_{-2} \zeta _2
   \nonumber \\
   &&
   -s_{-1} \zeta _3+\frac{s_1 \zeta
   _3}{2}+\frac{1}{2} \zeta _2 \bar{s}_{-2}-\zeta _3 \bar{s}_{-1}+\frac{5}{8} \zeta _3
   \bar{s}_1+\frac{1}{2} \zeta _2 \bar{s}_2-\zeta _2 s_{1,-1}+s_{3,-1}-\bar{s}_{-3,1}
   \nonumber \\
   &&
   -\zeta _2
   \bar{s}_{1,-1}-\zeta _2 \bar{s}_{1,1}-2 s_{1,2,-1}-2 s_{2,1,-1}+\bar{s}_{-2,1,1}+2
   \bar{s}_{1,-2,1}+\bar{s}_{2,-1,1}
   \nonumber \\
   &&+3 s_{1,1,1,-1}
   -\bar{s}_{1,-1,1,1}-2
   \bar{s}_{1,1,-1,1}
   \end{eqnarray}

     \begin{eqnarray}
     s_{1,-1} \bar{s}_{-1,-1} & = &  -\frac{5 \ln _2^4}{24}-\frac{2}{3}
   s_{-1} \ln _2^3+\frac{2}{3} s_1 \ln _2^3-\frac{2}{3} \bar{s}_{-1} \ln _2^3+\frac{2}{3}
   \bar{s}_1 \ln _2^3-\frac{1}{2} s_{-2} \ln _2^2+\frac{1}{2} s_2 \ln _2^2
   \nonumber \\
   &&
   -\frac{9}{4} \zeta _2
   \ln _2^2-\frac{1}{2} s_{-1,-1} \ln _2^2+\frac{1}{2} s_{-1,1} \ln _2^2+2 s_{1,-1} \ln
   _2^2+\frac{1}{2} \bar{s}_{-1,-1} \ln _2^2-\frac{3}{2} \bar{s}_{-1,1} \ln _2^2
   \nonumber \\
   &&
   -s_{-3} \ln _2+s_3
   \ln _2-s_{-1} \zeta _2 \ln _2
   -\frac{3}{2} s_1 \zeta _2 \ln _2+\bar{s}_{-3} \ln _2+\frac{1}{2}
   \zeta _2 \bar{s}_{-1} \ln _2
   \nonumber \\
   &&
   +\frac{1}{2} \zeta _2 \bar{s}_1 \ln _2-\bar{s}_3 \ln _2-2 s_{-2,-1}
   \ln _2-s_{-1,-2} \ln _2+s_{-1,2} \ln _2+s_{1,-2} \ln _2
    \nonumber \\
   &&
   -s_{1,2} \ln _2+2 \bar{s}_{-2,-1} \ln
   _2+\bar{s}_{-1,-2} \ln _2-\bar{s}_{-1,2} \ln _2-\bar{s}_{1,-2} \ln _2
   +\bar{s}_{1,2} \ln _2
    \nonumber \\
   &&   
   +2
   s_{-1,1,-1} \ln _2+2 s_{1,-1,-1} \ln _2-2 \bar{s}_{-1,1,-1} \ln _2-2 \bar{s}_{1,-1,-1} \ln
   _2-\frac{\zeta _2^2}{5}
    \nonumber \\
   &&
   +\text{LiHalf}_4+\frac{1}{2} s_{-2} \zeta _2-\frac{1}{8} s_{-1} \zeta
   _3+\frac{s_1 \zeta _3}{2}-\frac{1}{2} \zeta _2 \bar{s}_{-2}-\frac{1}{2} \zeta _3
   \bar{s}_{-1}-\frac{3}{4} \zeta _3 \bar{s}_1
    \nonumber \\
   &&
   -\frac{1}{2} \zeta _2 s_{-1,1}-\frac{1}{2} \zeta _2
   s_{1,-1}+s_{3,-1}+\frac{1}{2} \zeta _2 \bar{s}_{-1,1}+\frac{1}{2} \zeta _2
   \bar{s}_{1,-1}-\bar{s}_{3,-1}
    \nonumber \\
   &&
   -s_{-2,-1,-1}-s_{-1,-2,-1}-s_{1,2,-1}-s_{2,1,-1}+2
   \bar{s}_{-2,-1,-1}+\bar{s}_{-1,-2,-1}
    \nonumber \\
   &&
   +\bar{s}_{1,2,-1}+s_{-1,-1,1,-1}+s_{-1,1,-1,-1}+s_{1,-1,-1
   ,-1}
    \nonumber \\
   &&
   -\bar{s}_{-1,1,-1,-1}-2 \bar{s}_{1,-1,-1,-1}
   \end{eqnarray}

    \begin{eqnarray}
    s_{1,1} \bar{s}_{1,1} & = &  \frac{12
   \zeta _2^2}{5}-s_2 \zeta _2-\bar{s}_2 \zeta _2+2 s_{1,1} \zeta _2+2 \bar{s}_{1,1} \zeta _2+3
   s_1 \zeta _3+3 \zeta _3 \bar{s}_1+s_{3,1}
    \nonumber \\
   &&
   +\bar{s}_{3,1}-2 s_{1,2,1}-2 s_{2,1,1}-2
   \bar{s}_{1,2,1}-2 \bar{s}_{2,1,1}+3 s_{1,1,1,1}+3 \bar{s}_{1,1,1,1}
   \end{eqnarray}

     \begin{eqnarray}
     s_{1,1}
   \bar{s}_{-1,-1} & = &  -\frac{\ln _2^4}{3}-\frac{2}{3} s_{-1} \ln _2^3+\frac{2}{3} s_1 \ln
   _2^3-\frac{1}{3} \bar{s}_{-1} \ln _2^3+\frac{1}{3} \bar{s}_1 \ln _2^3-\frac{1}{2} s_{-2} \ln
   _2^2+\frac{1}{2} s_2 \ln _2^2
    \nonumber \\
   &&
   +\zeta _2 \ln _2^2-\frac{1}{2} \bar{s}_{-2} \ln _2^2+\frac{1}{2}
   \bar{s}_2 \ln _2^2-\frac{1}{2} s_{-1,-1} \ln _2^2+\frac{1}{2} s_{-1,1} \ln _2^2+\frac{3}{2}
   s_{1,-1} \ln _2^2
    \nonumber \\
   &&
   +\frac{1}{2} s_{1,1} \ln _2^2+\frac{1}{2} \bar{s}_{-1,-1} \ln _2^2-\frac{1}{2}
   \bar{s}_{-1,1} \ln _2^2+\frac{1}{2} \bar{s}_{1,-1} \ln _2^2-\frac{1}{2} \bar{s}_{1,1} \ln
   _2^2+s_{-1} \zeta _2 \ln _2
    \nonumber \\
   &&
   +s_1 \zeta _2 \ln _2+\frac{\zeta _3 \ln _2}{4}+\bar{s}_{-3} \ln
   _2-\bar{s}_3 \ln _2-s_{-2,-1} \ln _2-s_{-2,1} \ln _2+\bar{s}_{-2,-1} \ln _2
    \nonumber \\
   &&
   -\bar{s}_{-2,1} \ln
   _2
   +\bar{s}_{-1,-2} \ln _2-\bar{s}_{-1,2} \ln _2-\bar{s}_{1,-2} \ln _2+\bar{s}_{1,2} \ln
   _2+s_{-1,1,-1} \ln _2
    \nonumber \\
   &&
   +s_{-1,1,1} \ln _2+s_{1,-1,-1} \ln _2+s_{1,-1,1} \ln _2-\bar{s}_{-1,1,-1}
   \ln _2+\bar{s}_{-1,1,1} \ln _2-\bar{s}_{1,-1,-1} \ln _2
    \nonumber \\
   &&
   +\bar{s}_{1,-1,1} \ln _2
   +\frac{3 \zeta
   _2^2}{20}-\frac{1}{2} s_{-2} \zeta _2-\frac{s_2 \zeta _2}{2}+\frac{1}{8} s_{-1} \zeta
   _3+\frac{s_1 \zeta _3}{4}+\frac{1}{4} \zeta _3 \bar{s}_{-1}+\frac{1}{8} \zeta _3
   \bar{s}_1
    \nonumber \\
   &&
   +\frac{1}{2} \zeta _2 \bar{s}_2+\frac{1}{2} \zeta _2 s_{-1,-1}+\frac{1}{2} \zeta _2
   s_{-1,1}+\frac{1}{2} \zeta _2 s_{1,-1}+\frac{1}{2} \zeta _2
   s_{1,1}+s_{3,1}+\bar{s}_{-3,-1}
    \nonumber \\
   &&-\frac{1}{2} \zeta _2 \bar{s}_{-1,-1}-\frac{1}{2} \zeta _2
   \bar{s}_{1,1}-s_{-2,-1,1}-s_{-1,-2,1}-s_{1,2,1}-s_{2,1,1}-\bar{s}_{-2,1,-1}
    \nonumber \\
   &&
   -\bar{s}_{-1,2,-1}-
   \bar{s}_{1,-2,-1}-\bar{s}_{2,-1,-1}+s_{-1,-1,1,1}+s_{-1,1,-1,1}+s_{1,-1,-1,1}
    \nonumber \\
   &&
   +\bar{s}_{-1,1,1,-1
   }+\bar{s}_{1,-1,1,-1}+\bar{s}_{1,1,-1,-1}
   \end{eqnarray}

     \begin{eqnarray}
     s_{-1,-1} \bar{s}_{-1,-1} & = &  -\ln _2^4+2
   \zeta _2 \ln _2^2+2 s_{-1,-1} \ln _2^2+2 \bar{s}_{-1,-1} \ln _2^2+s_{-3} \ln _2-s_3 \ln
   _2
    \nonumber \\
   &&
   -s_{-1} \zeta _2 \ln _2-4 \zeta _3 \ln _2+\bar{s}_{-3} \ln _2-\zeta _2 \bar{s}_{-1} \ln
   _2-\bar{s}_3 \ln _2+2 s_{-1,-2} \ln _2
    \nonumber \\
   &&
   -2 s_{-1,2} \ln _2
   -2 s_{2,-1} \ln _2+2 \bar{s}_{-1,-2}
   \ln _2-2 \bar{s}_{-1,2} \ln _2-2 \bar{s}_{2,-1} \ln _2
    \nonumber \\
   &&
   +4 s_{-1,-1,-1} \ln _2+4
   \bar{s}_{-1,-1,-1} \ln _2+\frac{14 \zeta _2^2}{5}-8 \text{LiHalf}_4+\frac{s_2 \zeta
   _2}{2}
    \nonumber \\
   &&
   +\frac{5}{4} s_{-1} \zeta _3+\frac{5}{4} \zeta _3 \bar{s}_{-1}+\frac{1}{2} \zeta _2
   \bar{s}_2+s_{-3,-1}-\zeta _2 s_{-1,-1}+\bar{s}_{-3,-1}
    \nonumber \\
   &&
   -\zeta _2 \bar{s}_{-1,-1}-2 s_{-1,2,-1}-2
   s_{2,-1,-1}-2 \bar{s}_{-1,2,-1}-2 \bar{s}_{2,-1,-1}
    \nonumber \\
   &&
   +3 s_{-1,-1,-1,-1}+3
   \bar{s}_{-1,-1,-1,-1} 
   \end{eqnarray}


%
%

%



\begin{thebibliography}{}
%
%
%

\bibitem{Prygarin:2018tng} 
  A.~Prygarin,
  ``Reflection identities of harmonic sums up to weight three,''
  arXiv:1808.09307 [hep-th].







\bibitem{handbuch}
 N.~Nielsen, Handbuch der Theorie der Gammafunktion, (Teubner, Leipzig, 1906); reprinted by Chelsea Publishing Company, Bronx, NY, 1965.



\bibitem{Blumlein:2009ta} 
  J.~Blumlein,
  ``Structural Relations of Harmonic Sums and Mellin Transforms up to Weight $w = 5$,''
  Comput.\ Phys.\ Commun.\  {\bf 180}, 2218 (2009)
  doi:10.1016/j.cpc.2009.07.004
  [arXiv:0901.3106 [hep-ph]].
  
  
  

\bibitem{HS1}
A. Gonzalez-Arroyo, C. Lopez, and F.J. Yndurain,
Nucl. Phys. {\bf B153} (1979) 161;\\
A. Gonzalez-Arroyo and C. Lopez,
Nucl. Phys. {\bf B166} (1980) 429.




\bibitem{Vermaseren:1998uu} 
  J.~A.~M.~Vermaseren,
  ``Harmonic sums, Mellin transforms and integrals,''
  Int.\ J.\ Mod.\ Phys.\ A {\bf 14}, 2037 (1999)
  doi:10.1142/S0217751X99001032
  [hep-ph/9806280].
  
  
  

\bibitem{Blumlein:1998if} 
  J.~Blumlein and S.~Kurth,
  ``Harmonic sums and Mellin transforms up to two loop order,''
  Phys.\ Rev.\ D {\bf 60}, 014018 (1999)
  doi:10.1103/PhysRevD.60.014018
  [hep-ph/9810241].
  
  
    
\bibitem{Remiddi:1999ew} 
  E.~Remiddi and J.~A.~M.~Vermaseren,
  ``Harmonic polylogarithms,''
  Int.\ J.\ Mod.\ Phys.\ A {\bf 15}, 725 (2000)
  doi:10.1142/S0217751X00000367
  [hep-ph/9905237].
  
  
  
\bibitem{Kotikov:2005gr} 
  A.~V.~Kotikov and V.~N.~Velizhanin,
  ``Analytic continuation of the Mellin moments of deep inelastic structure functions,''
  hep-ph/0501274.
  
  
  
 
  
  
  

\bibitem{Gromov:2015vua} 
  N.~Gromov, F.~Levkovich-Maslyuk and G.~Sizov,
  ``Pomeron Eigenvalue at Three Loops in $\mathcal N=$ 4 Supersymmetric Yang-Mills Theory,''
  Phys.\ Rev.\ Lett.\  {\bf 115}, no. 25, 251601 (2015)
  doi:10.1103/PhysRevLett.115.251601
  [arXiv:1507.04010 [hep-th]].
  


\bibitem{Alfimov:2018cms} 
  M.~Alfimov, N.~Gromov and G.~Sizov,
  ``BFKL spectrum of $ \mathcal{N} $ = 4: non-zero conformal spin,''
  JHEP {\bf 1807}, 181 (2018)
  doi:10.1007/JHEP07(2018)181
  [arXiv:1802.06908 [hep-th]].
  
    
\bibitem{Caron-Huot:2016tzz} 
  S.~Caron-Huot and M.~Herranen,
  ``High-energy evolution to three loops,''
  JHEP {\bf 1802}, 058 (2018)
  doi:10.1007/JHEP02(2018)058
  [arXiv:1604.07417 [hep-ph]].
   
\bibitem{Kotikov:2006ts} 
  A.~V.~Kotikov and L.~N.~Lipatov,
  ``On the highest transcendentality in N=4 SUSY,''
  Nucl.\ Phys.\ B {\bf 769}, 217 (2007)
  doi:10.1016/j.nuclphysb.2007.01.020
  [hep-th/0611204].
  

  
  
  
   
  
  


\bibitem{AblingerThesis}
J.~Ablinger, A Computer Algebra Toolbox for Harmonic Sums Related to Particle Physics, Diploma Thesis, arXiv:1011.1176 [math-ph];\\
J.~Ablinger, Computer Algebra Algorithms for Special Functions in Particle Physics, PhD-Thesis, Johannes Kepler University Linz, April 2012


\bibitem{Ablinger:2011te} 
  J.~Ablinger, J.~Blumlein and C.~Schneider,
  ``Harmonic Sums and Polylogarithms Generated by Cyclotomic Polynomials,''
  J.\ Math.\ Phys.\  {\bf 52}, 102301 (2011)
  doi:10.1063/1.3629472
  [arXiv:1105.6063 [math-ph]].


  
  
\bibitem{Blumlein:2009fz} 
  J.~Blümlein,
  ``Structural Relations of Harmonic Sums and Mellin Transforms at Weight $w= 6$,''
  Clay Math.\ Proc.\  {\bf 12}, 167 (2010)
  [arXiv:0901.0837 [math-ph]].
  
  
  




  


































  
\bibitem{Bondarenko:2015tba} 
  S.~Bondarenko and A.~Prygarin,
  ``Hermitian separability and transition from singlet to adjoint BFKL equations in $\mathcal{N}=4$ super Yang-Mills Theory,''
  arXiv:1510.00589 [hep-th].
 
\bibitem{Bondarenko:2016tws} 
  S.~Bondarenko and A.~Prygarin,
  ``On a residual freedom of the next-to-leading BFKL eigenvalue in color adjoint representation in planar $ \mathcal{N}=4 $ SYM,''
  JHEP {\bf 1607}, 081 (2016)
  doi:10.1007/JHEP07(2016)081
  [arXiv:1603.01093 [hep-th]].
  























 


\bibitem{Maitre:2005uu} 
  D.~Maitre,
  ``HPL, a mathematica implementation of the harmonic polylogarithms,''
  Comput.\ Phys.\ Commun.\  {\bf 174}, 222 (2006)
  doi:10.1016/j.cpc.2005.10.008
  [hep-ph/0507152].
  
  
  
  
  
  



  
  
\end{thebibliography}

\end{document}